\documentclass[11pt]{amsart}
\usepackage{amssymb,amsmath,a4wide}
\usepackage[pagebackref=true]{hyperref}
\usepackage[all]{xy}

\title[Homotopical interpretation of globular complex]{Homotopical
 interpretation of globular complex by multipointed d-space}

\author[P. Gaucher]{Philippe Gaucher}
\address{Laboratoire PPS  (CNRS UMR 7126)\\ Universit{\'e} Paris 7--Denis Diderot\\
  Case 7014\\ 75205 PARIS Cedex 13 \\ France}
\email{gaucher@pps.jussieu.fr}
\urladdr{http://www.pps.jussieu.fr/{\~{}}gaucher/}
\subjclass{55U35, 18G55, 55P99, 68Q85}
\keywords{homotopy, directed homotopy, combinatorial model category,
  simplicial category, topological category, delta-generated
  space, d-space, globular complex, time flow}

\swapnumbers

\newcommand{\C}{\mathcal{C}}

\newcommand{\Z}{\mathbb{Z}}

\newcommand{\R}{\mathbb{R}}

\newcommand{\p}{\times}
\renewcommand{\vec}{\overrightarrow}
\renewcommand{\P}{\mathbb{P}}

\newcommand{\be}{\begin{equation}}
\newcommand{\ee}{\end{equation}}
\newcommand{\bea}{\begin{eqnarray}}
\newcommand{\eea}{\end{eqnarray}}
\newcommand{\beas}{\begin{eqnarray*}}
\newcommand{\eeas}{\end{eqnarray*}}

\newtheorem{thm}{Theorem}[section]
\newtheorem{prop}[thm]{Proposition}

\newtheorem{cor}[thm]{Corollary}

\newtheorem*{conv}{Convention}
\newtheorem{conj}[thm]{Conjecture}
\newtheorem{defn}[thm]{Definition}
\newtheorem{nota}[thm]{Notation}

\newcommand{\bd}{\begin{defn}}
\newcommand{\ed}{\end{defn}}
\newcommand{\bp}{\begin{prop}}
\newcommand{\ep}{\end{prop}}
\newcommand{\bth}{\begin{thm}}
\renewcommand{\eth}{\end{thm}}
\newcommand{\bpf}{\begin{proof}}
\newcommand{\epf}{\end{proof}}

\newcommand{\fl}[1]{\ar@{->}[l]_-{#1}}
\newcommand{\fr}[1]{\ar@{->}[r]^-{#1}}
\newcommand{\fd}[1]{\ar@{->}[d]_-{#1}}
\newcommand{\fu}[1]{\ar@{->}[u]^-{#1}}

\renewcommand{\top}{{\mathbf{Top}}}
\newcommand{\cgtop}{{\mathbf{CGTop}}}
\newcommand{\iso}{\cong}

\newcommand{\ot}{\otimes}
\DeclareMathOperator{\sing}{Sing}

\newcommand{\vI}{\vec{I}}
\renewcommand{\leq}{\leqslant}
\renewcommand{\geq}{\geqslant}

\def\cartesien{%
  \ar@{-}[]+R+<6pt,-2pt>;[]+RD+<6pt,-6pt>%
  \ar@{-}[]+D+<2pt,-6pt>;[]+RD+<6pt,-6pt>%
}
\def\cocartesien{%
  \ar@{-}[]+L+<-6pt,+2pt>;[]+LU+<-6pt,+6pt>%
  \ar@{-}[]+U+<-2pt,+6pt>;[]+LU+<-6pt,+6pt>%
}
\def\hocartesien{%
  \ar@{-}[]+R+<6pt,-2pt>;[]+RD+<6pt,-6pt>_{h}%
  \ar@{-}[]+D+<2pt,-6pt>;[]+RD+<6pt,-6pt>%
}
\def\hococartesien{%
  \ar@{-}[]+L+<-6pt,+2pt>;[]+LU+<-6pt,+6pt>_{h}%
  \ar@{-}[]+U+<-2pt,+6pt>;[]+LU+<-6pt,+6pt>%
}

\newcommand{\brm}[1]{\rm{\mathbf{#1}}} 
\newcommand{\ho}{{\brm{Ho}}}
\newcommand{\gltop}{{\brm{glTop}}}
\newcommand{\dtop}{{\brm{Flow}}}
\newcommand{\set}{{\brm{Set}}}

\newcommand{\ttop}{{\brm{TOP}}}

\newcommand{\mtop}{{\brm{MTop}}}
\newcommand{\mttop}{{\brm{MTOP}}}
\newcommand{\mdtop}{{\brm{MdTop}}}
\newcommand{\mdttop}{{\brm{MDTOP}}}
\newcommand{\glob}{{\rm{Glob}}}
\newcommand{\map}{{\rm{Map}}}
\newcommand{\res}{{\!}\restriction}

\newcommand{\liminj}{\varinjlim}
\newcommand{\limproj}{\varprojlim}

\makeatletter

\def\varholim@#1#2{%
  \vtop{\m@th\ialign{##\cr
    \hfil$#1\operator@font holim$\hfil\cr
    \noalign{\nointerlineskip\kern1.5\ex@}#2\cr
    \noalign{\nointerlineskip\kern-\ex@}\cr}}%
}
\def\holimproj{%
  \mathop{\mathpalette\varholim@{\leftarrowfill@\textstyle}}\nmlimits@
}
\def\holiminj{%
  \mathop{\mathpalette\varholim@{\rightarrowfill@\textstyle}}\nmlimits@
}

\makeatother

\DeclareMathOperator{\id}{Id}

\newcommand{\cell}{{\brm{cell}}}
\newcommand{\cof}{{\brm{cof}}}
\newcommand{\inj}{{\brm{inj}}}
\DeclareMathOperator{\Mod}{Mod}

\begin{document}

\begin{abstract} 
  Globular complexes were introduced by E. Goubault and the author to
  model higher dimensional automata.  Globular complexes are
  topological spaces equipped with a globular decomposition which is
  the directed analogue of the cellular decomposition of a CW-complex.
  We prove that there exists a combinatorial model category such that
  the cellular objects are exactly the globular complexes and such
  that the homotopy category is equivalent to the homotopy category of
  flows.  The underlying category of this model category is a variant
  of M. Grandis' notion of d-space over a topological space colimit
  generated by simplices. This result enables us to understand the
  relationship between the framework of flows and other works in
  directed algebraic topology using d-spaces. It also enables us to
  prove that the underlying homotopy type functor of flows can be
  interpreted up to equivalences of categories as the total left
  derived functor of a left Quillen adjoint.
\end{abstract}

\maketitle

\tableofcontents

\section{Introduction}

Globular complexes were introduced by E. Goubault and the author to
model higher dimensional automata in \cite{diCW}, and studied further
in \cite{model2}. They are topological spaces modeling a state space
equipped with a globular decomposition encoding the temporal ordering
which is a directed analogue of the cellular decomposition of a
CW-complex.

The fundamental geometric shape of this topological model of
concurrency is the topological globe of a space $Z$, practically of a
$n$-dimensional disk for some $n\geq 0$. The underlying state space is
the quotient
\[{\frac{\{\widehat{0},\widehat{1}\}\sqcup (Z\p[0,1])}{(z,0) = (z',0)
    = \widehat{0},(z,1) = (z',1)=\widehat{1}}}\] equal to the
unreduced suspension of $Z$ if $Z\neq \varnothing$ and equal to the
discrete space $\{\widehat{0},\widehat{1}\}$ if the space $Z$ is
empty. The segment $[0,1]$ together with the usual total ordering
plays in this setting the role of time ordering. The point
$\widehat{0}$ is the initial state and the point $\widehat{1}$ is the
final state of the globe of $Z$.  The execution paths are the
continuous maps $t\mapsto (z,t)$ and all strictly increasing
reparametrizations preserving the initial and final states. This
construction means that, between $\widehat{0}$ and $\widehat{1}$,
orthogonally to the time flow, there is, up to strictly increasing
reparametrization, a topological space of execution paths $Z$ which
represents the geometry of concurrency.  By pasting together this kind
of geometric shape with $Z$ being a sphere or a disk using attaching
maps locally preserving time ordering, it is then possible to
construct, up to homotopy, any time flow of any concurrent process.
In particular, the time flow of any process algebra can be modeled by
a precubical set \cite{ccsprecub}, and then by a globular complex
using a realization functor from precubical sets to globular complexes
(\cite[Proposition~3.9]{diCW} and \cite[Theorem~5.4.2]{realization}).
See also in \cite{diCW} several examples of PV diagrams whose first
appearance in computer science goes back to \cite{EWDCooperating}, and
in concurrency theory to \cite{Gunawardena1}.

Although the topological model of globular complexes is included in
all other topological models \cite{survol} introduced for this purpose
(local pospace \cite{MR1683333}, $d$-space \cite{mg}, stream
\cite{SK}), it is therefore expressive enough to contain all known
examples coming from concurrency.

However, the category of globular complexes alone does not satisfy any
good mathematical property for doing homotopy because it is, in a
sense, too small. In particular, it is not complete nor cocomplete.
This is one of the reasons for introducing the category of flows in
\cite{model3} and for constructing a functor associating a globular
complex with a flow in \cite{model2} allowing the interpretation of
some geometric properties of globular complexes in the model category
of flows.

We prove in this work that a variant of M. Grandis' notion of
$d$-space \cite{mg} can be used to give a homotopical interpretation
of the notion of globular complex.  Indeed, using this variant, it is
possible to construct a combinatorial model category such that the
globular complexes are exactly the cellular objects.

This result must be understood as a directed version of the following
fact: The category of cellular spaces, in which the cells are not
necessarily attached by following the increasing ordering of
dimensions as in CW-complexes, is the category of cellular objects of
the usual model category of topological spaces.  Moreover, if we
choose to work in the category of $\Delta$-generated spaces, i.e.
with spaces which are colimits of simplices, then the model category
becomes combinatorial. 

It turns out that the model category of multipointed $d$-spaces has a
homotopy category which is equivalent to the homotopy category of
flows of \cite{model3}. So the result of this paper enables us to
understand the relationship between the framework of flows and other
works in directed algebraic topology using M. Grandis' $d$-spaces.

As a straightforward application, it is also proved that, up to
equivalences of categories, the underlying homotopy type functor of
flows introduced in \cite{model2} can be viewed as the total left
derived functor of a left Quillen adjoint. This result is interesting
since this functor is complicated to use.  Indeed, it takes a flow to
a homotopy type of topological space. The latter plays the role of the
underlying state space which is unique only up to homotopy, not up to
homeomorphism, in the framework of flows. This result will simplify
future calculations of the underlying homotopy type thanks to the
possibility of using homotopy colimit techniques.

\subsection*{Outline of the paper}

Section~\ref{topp} is devoted to a short exposition about topological
spaces ($k$-space, $\Delta$-generated space, compactly generated
space). Proposition~\ref{connexe} and Proposition~\ref{bientop} seem
to be new.  Section~\ref{multi} presents the variant of Grandis'
notion of $d$-space which is used in the paper and it is proved that
this new category is locally presentable.  The new, and short,
definition of a globular complex is given in Section~\ref{new}. It is
also proved in this section that the new globular complexes are
exactly the ones previously defined in \cite{model2}.
Section~\ref{retract} is a technical section which sketches the theory
of inclusions of a strong deformation retract in the category of
multipointed $d$-spaces. The main result, the closure of these maps
under pushout, is used in the construction of the model structure.
Section~\ref{comb} constructs the model structure.  Section~\ref{comp}
establishes the equivalence between the homotopy category of
multipointed $d$-spaces and the homotopy category of flows of
\cite{model3}. The same section also explores other connections
between multipointed $d$-spaces and flows.  In particular, it is
proved that there is a kind of left Quillen equivalence of cofibration
categories from multipointed $d$-spaces to flows. And finally
Section~\ref{under} is the application interpreting the underlying
homotopy type functor of flows as a total left derived functor.

\subsection*{Prerequisites}

There are many available references for general topology applied to
algebraic topology, e.g., \cite{topologie} \cite{MR1867354}. However,
the notion of $k$-space which is presented is not exactly the good
one. In general, the category of $k$-spaces is unfortunately defined
as the coreflective hull of the category of quasi-compact spaces (i.e.
spaces satisfying the finite open subcovering property and which are
not necessarily Hausdorff) which is not cartesian closed
(\cite{breger}, and \cite[Theorem~3.6]{MR1135098}).  One then obtains
a cartesian closed full subcategory by restricting to Hausdorff spaces
(\cite[Definition~7.2.5 and Corollary~7.2.6]{MR96g:18001b}).  However,
it is preferable to use the notion of weak Hausdorff space since some
natural constructions can lead outside this category.  So
\cite[Chapter~5]{concise} or \cite{MR1074175} Appendix~A must be
preferred for a first reading. See also \cite{MR2273730} and the
appendix of \cite{Ref_wH}.  Section~\ref{topp} of this paper is an
important section collecting the properties of topological spaces used
in this work. In particular, the category of $k$-spaces is defined as
the coreflective hull of the full subcategory of compact spaces, i.e.
of quasi-compact Hausdorff spaces. The latter category is cartesian
closed.

The reading of this work requires some familiarity with model category
techniques \cite{MR99h:55031} \cite{ref_model2}, with category theory
\cite{MR1712872} \cite{MR96g:18001a}\cite{gz}, and especially with
locally presentable categories \cite{MR95j:18001} and topological
categories \cite{topologicalcat}. 

\subsection*{Notations}

All categories are locally small. The set of morphisms from $X$ to $Y$
in a category $\C$ is denoted by $\C(X,Y)$. The identity of $X$ is
denoted by $\id_X$. Colimits are denoted by $\liminj$ and limits by
$\limproj$. Let $\C$ be a cocomplete category. The class of morphisms
of $\C$ that are transfinite compositions of pushouts of elements of a
set of morphisms $K$ is denoted by $\cell(K)$. An element of
$\cell(K)$ is called a \textit{relative $K$-cell complex}.  The
category of sets is denoted by $\set$. The class of maps satisfying
the right lifting property with respect to the maps of $K$ is denoted
by $\inj(K)$. The class of maps satisfying the left lifting property
with respect to the maps of $\inj(K)$ is denoted by $\cof(K)$. The
cofibrant replacement functor of a model category is denoted by
$(-)^{cof}$. The notation $\simeq$ means \textit{weak equivalence} or
\textit{equivalence of categories}, the notation $\iso$ means
\textit{isomorphism}. A combinatorial model category is a cofibrantly
generated model category such that the underlying category is locally
presentable. The notation $\id_A$ means identity of $A$.  The initial
object (resp. final object) of a category is denoted by $\varnothing$
(resp.  $\mathbf{1}$). In a cofibrantly generated model category with
set of generating cofibrations $I$, a \textit{cellular object} is an
object $X$ such that the map $\varnothing \rightarrow X$ belongs to
$\cell(I)$. The cofibrant objects are the retracts of the cellular
objects in a cofibrantly generated model category.

\subsection*{Acknowledgments} I thank very much Ji\v{r}{\'{\i}} Rosick\'y
for answering my questions about topological and locally presentable
categories.

\section{About topological spaces}
\label{topp}

We must be very careful in this paper since we are going to work with
$\Delta$-generated spaces without any kind of separation condition.
However, every compact space is Hausdorff. 

Let $\mathcal{T}$ be the category of general topological spaces
\cite{topologie} \cite{MR1867354}. This category is complete and
cocomplete. Limits are obtained by taking the initial topology, and
colimits are obtained by taking the final topology on the underlying
(co)limits of sets. This category is the paradigm of topological
category because of the existence of the initial and final structures
\cite{topologicalcat}.

A one-to-one continuous map $f:X\rightarrow Y$ between general
topological spaces is a \textit{(resp. closed) inclusion of spaces} if
$f$ induces a homeomorphism $X\iso f(X)$ where $f(X)$ is a (resp.
closed) subset of $Y$ equipped with the relative topology.  If $f$ is
a closed inclusion and if moreover for all $y\in Y\backslash f(X)$,
$\{y\}$ is closed in $Y$, then $f:X\rightarrow Y$ is called a
\textit{closed $T_1$-inclusion of spaces}.

\bp \label{ir} (well-known) Let $i:X\rightarrow Y$ be a continuous map
between general topological spaces. If there exists a retract
$r:Y\rightarrow X$, i.e. a continuous map $r$ with $r\circ i=\id_X$,
then $i$ is an inclusion of spaces and $X$ is equipped with the final
topology with respect to $r$, i.e. $r$ is a quotient map. \ep

\bpf The set map $i$ is one-to-one and the set map $r$ is onto.  Since
$i$ is the equalizer (resp. $p$ is the coequalizer) of the pair of
maps $(i\circ r,\id_Y)$, one has $X\iso \{y\in Y,y = i(r(y))\}\iso
i(X)$ and $X$ is a quotient of $Y$ equipped with the final topology.
\epf

Let us emphasize two facts related to Proposition~\ref{ir}:
\begin{enumerate}
\item There does not exist any reason for the map $i:X\rightarrow Y$
  to be a closed inclusion of spaces without any additional separation
  condition.
\item If the map $i:X\rightarrow Y$ is a closed inclusion of spaces
  anyway, then there does not exist any reason for the map
  $i:X\rightarrow Y$ to be a closed $T_1$-inclusion of spaces without
  any additional separation condition.
\end{enumerate}

Let $\mathcal{B}$ be a full subcategory of the category $\mathcal{T}$
of general topological spaces. A topological space $X$ is
\textit{$\mathcal{B}$-generated} if the natural map
\[k_\mathcal{B}(X):=\liminj_{B\in \mathcal{B},B\rightarrow X} B
\longrightarrow X\] is a homeomorphism (note that the diagram above
may be large). The underlying set of the space $k_\mathcal{B}(X)$ is
equal to the underlying set of the space $X$.  The
$\mathcal{B}$-generated spaces assemble to a full coreflective
subcategory of $\mathcal{T}$, in fact the coreflective hull of
$\mathcal{B}$, denoted by $\top_\mathcal{B}$.  The right adjoint to
the inclusion functor $\top_\mathcal{B} \subset \mathcal{T}$ is
precisely the functor $k_\mathcal{B}$, called the
\textit{Kelleyfication functor}.  The category $\top_\mathcal{B}$ is
complete and cocomplete.  Colimits in $\top_\mathcal{B}$ and in
$\mathcal{T}$ are the same.  Limits in $\top_\mathcal{B}$ are obtained
by calculating the limits in $\mathcal{T}$ and by applying the
Kelleyfication functor $k_\mathcal{B}$. See \cite{MR45:9323} for a
proof of all these facts.  The category $\top_\mathcal{B}$ is also
locally presentable as soon as $\mathcal{B}$ is small by
\cite[Theorem~3.6]{FR}. This fact was conjectured by J. H. Smith in
unpublished notes.

\begin{nota} The binary product in the category $\top_\mathcal{B}$ of
  $\mathcal{B}$-generated spaces is denoted by $\p_\mathcal{B}$, if
  necessary.  The binary product in the category $\mathcal{T}$ is
  denoted by $\p_\mathcal{T}$, if necessary. \end{nota}

\bp \label{inclusiontop} (\cite[Proposition~1.5]{MR45:9323}) Let
$\mathcal{B}_1$ and $\mathcal{B}_2$ be two full subcategories of the
category $\mathcal{T}$ with $\mathcal{B}_1 \subset \mathcal{B}_2$.
Then one has $\top_{\mathcal{B}_1} \subset \top_{\mathcal{B}_2}$ and
the inclusion functor $i_{\mathcal{B}_1}^{\mathcal{B}_2}$ has a right
adjoint $k_{\mathcal{B}_1}$.  \ep

\bd (e.g., \cite[Proposition~7.1.5]{MR96g:18001b} or \cite[p
160]{Ref_wH}) Let $X$ and $Y$ be two general topological spaces. The
space $\mathcal{T}_t(X,Y)$ is the set of continuous maps
$\mathcal{T}(X,Y)$ from $X$ to $Y$ equipped with the topology
generated by the subbasis $N(h,U)$ where $h:K\rightarrow X$ is a
continuous map from a compact $K$ to $X$, where $U$ is an open of $Y$
and where $N(h,U) = \{f:X\rightarrow Y,f(h(K))\subset U\}$. This
topology is called the {\rm compact-open topology}.  \ed

Proposition~\ref{cartesianclosed} is a slight modification of
\cite[Section~3]{MR45:9323}, which is stated without proof in
\cite{delta}. It is important because our theory requires a cartesian
closed category of topological spaces.

\bp (Dugger-Vogt) \cite{delta} \label{cartesianclosed} Let us suppose
that every object of $\mathcal{B}$ is compact and that the binary
product in $\mathcal{T}$ of two objects of $\mathcal{B}$ is
$\mathcal{B}$-generated. Then: 
\begin{enumerate}
\item For any $B\in \mathcal{B}$ and any $X\in \top_\mathcal{B}$, one
  has $B\p_\mathcal{T}X\in \top_\mathcal{B}$.
\item The category $\top_\mathcal{B}$ is cartesian closed.
\end{enumerate}
\ep

\bpf[Sketch of proof] We follow Vogt's proof.  It is well known that
for any compact space $K$, the functor $K\p_\mathcal{T}-:\mathcal{T}
\rightarrow \mathcal{T}$ has the right adjoint
$\mathcal{T}_t(K,-):\mathcal{T} \rightarrow \mathcal{T}$. So the
partial evaluation set map $\mathcal{T}_t(X,Y)\p_\mathcal{T} K
\rightarrow \mathcal{T}_t(K,Y)\p_\mathcal{T} K \rightarrow Y$ is
continuous for any general topological space $X$ and $Y$ and for any
continuous map $K\rightarrow X$ with $K$ compact.

Consider the topological space $\ttop_\mathcal{B}(X,Y) =
k_\mathcal{B}(\mathcal{T}_t(X,Y))$, i.e. the set
$\top_\mathcal{B}(X,Y)$ equipped with the Kelleyfication of the
compact-open topology. Let $X$ and $Y$ be two $\mathcal{B}$-generated
topological spaces. For any $B\in \mathcal{B}$, the composite
$g:B\rightarrow \mathcal{T}_t(X,Y)\p_\mathcal{T} X \rightarrow Y$ is
continuous for any continuous map $B\rightarrow
\mathcal{T}_t(X,Y)\p_\mathcal{T} X$ since $g$ is equal to the
composite $B\rightarrow \mathcal{T}_t(B,Y) \p_\mathcal{T} B
\rightarrow Y$ where $B\rightarrow\mathcal{T}_t(B,Y)$ is the composite
$B\rightarrow \mathcal{T}_t(X,Y)\p_\mathcal{T} X \rightarrow
\mathcal{T}_t(X,Y) \rightarrow \mathcal{T}_t(B,Y)$ and since $B$ is
compact. So the set map $\ttop_\mathcal{B}(X,Y)\p_\mathcal{B} X =
k_\mathcal{B}(\mathcal{T}_t(X,Y)\p_\mathcal{T} X) \rightarrow Y$ is
continuous.

By hypothesis, if $B_1,B_2\in \mathcal{B}$, then $B_1\p_\mathcal{T}
B_2\in \top_\mathcal{B}$. Since $B_1$ is compact, and since colimits
in $\top_\mathcal{B}$ and in $\mathcal{T}$ are the same, one deduces
that for any $X\in \top_\mathcal{B}$, $B_1\p_\mathcal{T} X\in
\top_\mathcal{B}$. So the canonical map $B_1\p_\mathcal{B} X =
k_\mathcal{B}(B_1\p_\mathcal{T} X) \rightarrow B_1\p_\mathcal{T} X$ is
a homeomorphism if $B_1\in \mathcal{B}$ and $X \in \top_\mathcal{B}$.
Hence the first assertion.

Let $X$, $Y$ and $Z$ be three $\mathcal{B}$-generated spaces. Let
$f:Y\p_\mathcal{B} X \rightarrow Z$ be a continuous map. Consider the
set map $\widetilde{f}: Y \rightarrow \ttop_\mathcal{B}(X,Z)$ defined
by $\widetilde{f}(y)(x) = f(y,x)$. Let $g:B\rightarrow Y$ be a
continuous map with $B\in \mathcal{B}$. Then the composite set map $B
\rightarrow Y \rightarrow \mathcal{T}_t(X,Z)$ is continuous since it
corresponds by adjunction to the continuous map $B\p_\mathcal{T}
X\stackrel{\iso} \rightarrow B\p_\mathcal{B} X \rightarrow
Y\p_\mathcal{B} X \rightarrow Z$.  Since $Y$ is
$\mathcal{B}$-generated, and therefore since $k_\mathcal{B}(Y)=Y$, the
set map $\widetilde{f}: Y \rightarrow \ttop_\mathcal{B}(X,Z)$ is
therefore continuous.  Hence the second assertion. \epf

For the rest of the section, let us suppose that $\mathcal{B}$
satisfies the following properties:
\begin{itemize} 
\item Every object of $\mathcal{B}$ is compact.
\item All simplices $\Delta^n=\{(t_0,\dots,t_n)\in
  (\R^{+})^n,t_0+\dots+t_n=1\}$ with $n\geq 0$ are objects of
  $\mathcal{B}$.
\item The binary product in $\mathcal{T}$ of two objects of
  $\mathcal{B}$ is $\mathcal{B}$-generated.
\end{itemize}

The category $\mathcal{K}$ of all compact spaces satisfies the
conditions above. It is the biggest possible choice. An object of
$\top_\mathcal{K}$ is called a \textit{$k$-space} \cite{concise}
\cite{MR1074175} \cite{Ref_wH}. 

The full subcategory $\Delta$ of $\mathcal{T}$ generated by all
topological simplices $\Delta^n = \{(t_0,\dots,t_n)\in
(\R^{+})^n,t_0+\dots+t_n=1\}$ with $n\geq 0$ is another possible
choice. Indeed, one has $\Delta^n\iso |\Delta[n]|$ where $|\Delta[n]|$
is the geometric realization of the $n$-simplex viewed as a simplicial
set. And there is a homeomorphism $\Delta^m \p_\mathcal{K} \Delta^n
\iso |\Delta[m]\p\Delta[n]|$ by \cite[Lemma~3.1.8]{MR99h:55031} or
\cite[Theorem~III.3.1 p49]{gz}~\footnote{Note that Gabriel and
  Zisman's proof is written in the full subcategory of Hausdorff
  spaces of the coreflective hull of the category of quasi-compact
  spaces, that is with the wrong notion of $k$-space.}.  At last,
since $\Delta^m$ is compact, the canonical map $\Delta^m
\p_\mathcal{K} \Delta^n \rightarrow \Delta^m \p_\mathcal{T} \Delta^n$
is an isomorphism. This choice is the smallest possible choice.
Further details about these topological spaces are available in
\cite{delta}.

As corollary of Proposition~\ref{cartesianclosed}, one has:

\bp \label{produit} Let $X$ and $Y$ be two $\mathcal{B}$-generated
spaces with $Y$ compact. Then the canonical map $X\p_\mathcal{B} Y
\rightarrow X\p_\mathcal{T} Y$ is a homeomorphism. \ep

\bpf Because the colimits in $\mathcal{T}$ and $\top_\mathcal{B}$ are
equal, one has:
\begin{align*}
  X\p_\mathcal{B} Y & \iso \liminj_{B\rightarrow X,B\in \mathcal{B}}
  (B\p_\mathcal{B} Y) & \hbox{ since $-\p_\mathcal{B} Y$ is a left
    adjoint by Proposition~\ref{cartesianclosed} (2)}\\
  & \iso \liminj_{B\rightarrow X,B\in \mathcal{B}} (B\p_\mathcal{T} Y)
  & \hbox{ since $B\p_\mathcal{T} Y \in \top_\mathcal{B}$ by
    Proposition~\ref{cartesianclosed} (1)}\\
  & \iso X\p_\mathcal{T} Y& \hbox{ since $Y$ is compact and since
    $-\p_\mathcal{T} Y$ is a left adjoint.}
\end{align*}
\epf 

The following two propositions exhibit two striking differences
between $\Delta$-generated spaces and $k$-spaces.

\bp \cite{Ref_wH} \cite{delta} \label{sous} Let $X$ be a general
topological space.
\begin{enumerate}
\item If $X$ is $\Delta$-generated, then every open subset of $X$
  equipped with the relative topology is $\Delta$-generated. There
  exists a closed subset of the topological $2$-simplex $\Delta^2$
  which is not $\Delta$-generated when it is equipped with the
  relative topology.
\item If $X$ is a $k$-space, then every closed subset of $X$ equipped
  with the relative topology is a $k$-space. An open subset of a
  $k$-space equipped with the relative topology need not be a
  $k$-space.
\end{enumerate}
\ep 

\bp \label{connexe} Every $\Delta$-generated space is homeomorphic to
the disjoint union of its non-empty connected components, which are
also its non-empty path-connected components.  In particular, a
$\Delta$-generated space is connected if and only if it is
path-connected. \ep

\bpf Let $X$ be a $\Delta$-generated space.  Let $\widehat{X}$ be the
disjoint sum of the connected components (resp. path-connected
components) of $X$ which is still a $\Delta$-generated space. There is
a canonical continuous map $\widehat{X} \rightarrow X$. Let $B\in
\Delta$ be a simplex. Then any continuous map $B \rightarrow X$
factors uniquely as a composite $B \rightarrow \widehat{X} \rightarrow
X$ since $B$ is connected (resp. path-connected). Thus $\widehat{X}
\iso k_\Delta(\widehat{X}) \iso k_\Delta(X) \iso X$.  \epf

For example, the set of rational numbers $\mathbb{Q}$ equipped with
the order topology is a non-discrete totally disconnected space. The
latter space is not $\Delta$-generated since $k_\Delta(\mathbb{Q})$ is
the set $\mathbb{Q}$ equipped with the discrete topology.

The category of $\mathcal{B}$-generated spaces contains among other
things all geometric realizations of simplicial sets, of cubical sets,
all simplicial and cubical complexes, the discrete spaces, the
$n$-dimensional sphere $\mathbf{S}^n$ and the $n$-dimensional disk
$\mathbf{D}^n$ for all $n\geq 0$, all their colimits in the category
of general topological spaces, and so all CW-complexes and all open
subspaces of these spaces equipped with the relative topology by
Proposition~\ref{sous} (1).

The category of spaces $\top_\mathcal{B}$ has other interesting
properties which will be useful for the paper.

\bp \label{bientop} The forgetful functor $\omega_\mathcal{B} :
\top_\mathcal{B} \rightarrow \set$ is topological and fibre-small in
the sense of \cite{topologicalcat}. \ep

\bpf The category $\top_\mathcal{B}$ is a concretely coreflective
subcategory of $\mathcal{T}$. Therefore it is topological by
\cite[Theorem~21.33]{topologicalcat}.  \epf

By the topological duality theorem (see
\cite[Theorem~21.9]{topologicalcat}), any cocone
$(f_i:\omega_\mathcal{B}(A_i) \rightarrow X)$ of $\top_\mathcal{B}$
has a unique $\omega_\mathcal{B}$-final lift $(\overline{f_i}:A_i
\rightarrow A)$. The space $A$ is the set $X$ equipped with the final
topology. In particular, as already mentioned, colimits in
$\top_\mathcal{B}$ and in $\mathcal{T}$ are the same.  Another
consequence is that any quotient of a $\mathcal{B}$-generated space
equipped with the final topology is also a $\mathcal{B}$-generated
space. In particular, by Proposition~\ref{ir}, any retract of a
$\mathcal{B}$-generated space is $\mathcal{B}$-generated.

\bd A general topological space $X$ is {\rm weak Hausdorff} if for any
continuous map $g:K\rightarrow X$ with $K$ compact, $g(K)$ is closed
in $X$. A {\rm compactly generated space} is a weak Hausdorff
$k$-space. \ed

Note that if $X$ is weak Hausdorff, then any subset of $X$ equipped
with the relative topology is weak Hausdorff and that
$k_\mathcal{B}(X)$ is weak Hausdorff as well since $k_\mathcal{B}(X)$
contains more closed subsets than $X$.

\begin{nota} The category of compactly generated topological spaces is
  denoted by $\cgtop$. \end{nota}

Let us conclude this section by some remarks explaining why the
homotopy theory of any of the preceding categories of spaces is the
same. By \cite{FR} \cite{MR99h:55031}, the inclusion functor
$i_\Delta^\mathcal{K}:\top_\Delta \rightarrow \top_\mathcal{K}$ is a
left Quillen equivalence for the Quillen model structure: it is
evident that it is a left Quillen adjoint; it is a Quillen equivalence
because the natural map $k_\Delta(X)\rightarrow X$ is a weak homotopy
equivalence since for any space $K\in \{\mathbf{S}^n,\mathbf{S}^n\p
[0,1],n\geq 0\}$, a map $K\rightarrow X$ factors uniquely as a
composite $K \rightarrow k_\Delta(X)\rightarrow X$.  The weak
Hausdorffization functor $H:\top_\mathcal{K} \rightarrow \cgtop$ left
adjoint to the inclusion functor $\cgtop\subset \top_\mathcal{K}$ is
also a left Quillen equivalence for the Quillen model structure by
\cite[Theorem~2.4.23 and Theorem~2.4.25]{MR99h:55031}. All cofibrant
spaces of $\top_\Delta$ and $\top_\mathcal{K}$ are weak Hausdorff, and
therefore compactly generated because of the inclusion $\top_\Delta
\subset \top_\mathcal{K}$.

\section{Multipointed $d$-spaces}
\label{multi}

\bd A {\rm multipointed space} is a pair $(|X|,X^0)$ where
\begin{itemize}
\item $|X|$ is a $\Delta$-generated space called the {\rm underlying
    space} of $X$.
\item $X^0$ is a subset of $|X|$ called the {\rm $0$-skeleton} of $X$.
\end{itemize}
A morphism of multipointed spaces $f:X=(|X|,X^0) \rightarrow
Y=(|Y|,Y^0)$ is a commutative square
\[
\xymatrix{
X^0 \fr{f^0}\fd{} & Y^0 \fd{} \\ 
|X| \fr{|f|} & |Y|.}
\] 
The corresponding category is denoted by $\mtop_\Delta$.  \ed

\begin{nota} \label{composition} Let $M$ be a topological space. Let
  $\phi_1$ and $\phi_2$ be two continuous maps from $[0,1]$ to $M$
  with $\phi_1(1)=\phi_2(0)$. Let us denote by $\phi_1 *_a \phi_2$
  (with $0<a<1$) the following continuous map: if $0\leq t\leq a$,
  $(\phi_1 *_a \phi_2)(t)=\phi_1(\frac{t}{a})$ and if $a\leq t\leq 1$,
  $(\phi_1 *_a \phi_2)(t)=\phi_2(\frac{t-a}{1-a})$.
\end{nota}

\bd A {\rm multipointed $d$-space $X$} is a triple
$(|X|,X^0,\P^{top}X)$ where
\begin{itemize}
\item $(|X|,X^0)$ is a multipointed space.
\item The set $\P^{top}X$ is the disjoint union of the sets
  $\P_{\alpha,\beta}^{top}X$ for $(\alpha,\beta)$ running over $X^0\p
  X^0$ where $\P_{\alpha,\beta}^{top}X$ is a set of continuous paths
  $\phi$ from $[0,1]$ to $|X|$ such that $\phi(0)=\alpha$ and
  $\phi(1)=\beta$ which is closed under composition $*_{1/2}$ and
  under strictly increasing reparametrization~\footnote{These two
    facts imply that the set of paths is closed under $*_t$ for any
    $t\in]0,1[$.}, that is for every $\phi\in
  \P^{top}_{\alpha,\beta}X$ and for any strictly increasing continuous
  map $\psi:[0,1]\rightarrow [0,1]$ with $\psi(0)=0$ and $\psi(1)=1$,
  $\phi\circ \psi\in \P^{top}_{\alpha,\beta}X$. The element of
  $\P^{top}X$ are called {\rm execution path} or {\rm $d$-path}.
\end{itemize}
A morphism of multipointed $d$-spaces $f:X=(|X|,X^0,\P^{top}X)
\rightarrow Y=(|Y|,Y^0,\P^{top}Y)$ is a commutative square
\[
\xymatrix{
X^0 \fr{f^0}\fd{} & Y^0 \fd{} \\ 
|X| \fr{|f|} & |Y|}
\] 
such that for every $\phi\in \P^{top}X$, one has $|f|\circ \phi\in
\P^{top}Y$.  Let $\P^{top}f(\phi) = |f|\circ \phi$. The corresponding
category is denoted by $\mdtop_\Delta$.  \ed

Any set $E$ can be viewed as one of the multipointed $d$-spaces
$(E,E,\varnothing)$, $(E,E,E)$ and $(E,\varnothing,\varnothing)$ with
$E$ equipped with the discrete topology.

\begin{conv} For the sequel, a set $E$ is always associated with the
  multipointed $d$-space $(E,E,\varnothing)$. \end{conv}

\bth \label{locpre} The category $\mdtop_\Delta$ is concrete
topological and fibre-small over the category of sets, i.e. the
functor $X\mapsto \omega_\Delta(|X|)$ from $\mdtop_\Delta$ to $\set$
is topological and fibre-small.  Moreover, the category
$\mdtop_\Delta$ is locally presentable.  \eth

\bpf We partially mimic the proof of \cite[Theorem~4.2]{FR} and we use
the terminology of \cite[Chapter~5]{MR95j:18001}. A multipointed
$d$-space $X$ is a concrete structure over the set
$\omega_\Delta(|X|)$ which is characterized by a $0$-skeleton, a
topology and a set of continuous paths satisfying several axioms. The
category $\top_\Delta$ is topological and fibre-small over the
category of sets by Proposition~\ref{bientop}. So by
\cite[Theorem~5.3]{MR629337}, it is isomorphic to the category of
models $\Mod(T)$ of a relational universal strict Horn theory $T$
without equality, i.e. all axioms are of the form $(\forall {x}),
\phi({x}) \Rightarrow \psi({x})$ where $\phi$ and $\psi$ are
conjunctions of atomic formulas without equalities. As in the proof of
\cite[Theorem~3.6]{FR}, let us suppose that $T$ contains all its
consequences. It means that a universal strict Horn sentence without
equality belongs to $T$ if and only if it holds for all models of
$T$. The theory $T$ may contain a proper class of axioms. Let us
construct as in the proof of \cite[Theorem~3.6]{FR} an increasing
chain of full coreflective subcategories such that the union is equal
to $\Mod(T) = \top_\Delta$
\[\Mod(T_0) \subset \Mod(T_1) \subset ... \Mod(T_\alpha) \subset
....\] where each $T_\alpha$ is a subset of axioms of $T$ indexed by
an ordinal $\alpha$. Let ${\alpha_0}$ be an ordinal such that the
subcategory $\Mod(T_{\alpha_0})$ contains all simplices $\Delta^n =
\{(t_0,\dots,t_n)\in (\R^{+})^n,t_0+\dots+t_n=1\}$ with $n\geq
0$. Then one has $\Mod(T) \subset \Mod(T_{\alpha_0})$ since the
category $\Mod(T_{\alpha_0})$ is cocomplete. So one obtains the
isomorphism of categories $\Mod(T_{\alpha_0})\iso \top_\Delta$. The
theory $T_{\alpha_0}$ is a universal strict Horn theory without
equality containing a set of axioms using $2^{\aleph_0}$-ary
relational symbols $R_j$ for $j\in J$ for some set $J$. See
\cite{MR588210} for a description of this theory for the category of
general topological spaces. One has the isomorphism of categories
$\mtop_\Delta \iso \Mod(T_{\alpha_0}\cup\{S\})$ where $S$ is a $1$-ary
relational symbol whose interpretation is the $0$-skeleton. Let us add
to the theory $T_{\alpha_0}\cup\{S\}$ a new $2^{\aleph_0}$-ary
relational symbol $R$ whose interpretation is the set of execution
paths. And let us add the axioms (as in the proof of
\cite[Theorem~4.2]{FR}):
\begin{enumerate}
\item $(\forall x) R(x) \Rightarrow (S(x_0) \wedge S(x_1))$
\item $(\forall x,y,z) \left(\left(\bigwedge_{0\leq t\leq 1/2} x_{2t}=z_t\right) \wedge
  \left(\bigwedge_{1/2\leq t\leq 1} y_{2t-1}=z_{\frac{1}{2}+t}\right) \wedge R(x)
  \wedge R(y)\right) \Rightarrow R(z)$
\item $(\forall x) R(x) \Rightarrow R(xt)$ where $t$ is a strictly
  increasing reparametrization
\item $(\forall x) R(x) \Rightarrow R_j(xa)$ where $j\in J$ and
  $T_{\alpha_0}$ satisfies $R_j$ for $a$.
\end{enumerate}
The first axiom says that for all execution paths $x$ of a
multipointed $d$-space $X$, one has $x_0,x_1\in X^0$. The second axiom
says that the set of execution paths is closed under composition
$*_{1/2}$, and the third one that the same set is closed under
strictly increasing reparametrization. The last axiom says that all
execution paths are continuous. Then the new theory $T'$ is a
relational universal strict Horn theory without equality containing a
set of axioms. So by \cite[Theorem~5.30]{MR95j:18001}, the category
$\mdtop_\Delta\iso \Mod(T')$ is locally presentable~\footnote{A
  relational universal strict Horn theory is a limit theory since the
  axiom $(\forall x)(A(x) \Rightarrow B(x))$ is equivalent to the
  axiom $(\forall x)(\exists y)(A(x) \Rightarrow (B(x)\wedge y=x))$:
  see \cite[p 324]{MR629337}.}, and by \cite[Theorem~5.3]{MR629337},
it is fibre-small concrete and topological over the category of
sets. \epf

Note that the same proof shows that the category $\mtop_\Delta$ is
also concrete topological and locally presentable. It is of course
possible to prove directly that $\mtop_\Delta$ and $\mdtop_\Delta$ are
concrete fibre-small and topological using Proposition~\ref{bientop}
and Proposition~\ref{limcolim} below. But we do not know how to prove
that they are locally presentable without using a logical argument.

Since all relational symbols appearing in the relational universal
strict Horn theory axiomatizing $\mdtop_\Delta$ are of arity at most
$2^{\aleph_0}$, the category $\mdtop_\Delta$ is locally
$\lambda$-presentable, where $\lambda$ denotes a regular cardinal
greater or equal to $2^{\aleph_0}$ (\cite[p 160]{MR1697766}).

The proof of the following proposition makes explicit, using the $*_t$
composition laws, the construction of colimits in $\mdtop_\Delta$.

\bp\label{limcolim} The functor $U_\Delta : X = (|X|,X^0,\P^{top}X)
\mapsto |X|$ from $\mdtop_\Delta$ to $\top_\Delta$ is topological and
fibre-small. In particular, it creates limits and colimits. \ep

\bpf That it is fibre-small is clear. It is very easy to prove that
the functor $U_\Delta$ is topological by observing that the functor
$\omega_\Delta\circ U_\Delta$ is topological by Theorem~\ref{locpre}
and that the functor $\omega_\Delta$ is topological by
Proposition~\ref{bientop}. We prefer to give a proof with an explicit
construction of the final lift using the $*_t$ composition laws
because this proof will be reused several times in the paper. By
\cite[Theorem~21.9]{topologicalcat}, it then suffices to see that the
forgetful functor $U_\Delta:\mdtop_\Delta \rightarrow \top_\Delta$
satisfies: for any cocone $(f_i:U_\Delta(A_i)\rightarrow X)$, there is
a unique $U_\Delta$-final lift $(\overline{f_i}:A_i \rightarrow A)$.
Let $A^0$ be the image of the set map $\liminj A_i^0 \rightarrow
X=|A|$.  Let $\P^{top}A$ be the set of all possible finite
compositions $*_t$ of $f_i\circ \phi_i$ with $\phi_i \in \P^{top}A_i$
and $t\in ]0,1[$.  Then one obtains a set of continuous paths which is
closed under strictly increasing reparametrization.  Indeed, let
$\phi_0 *_{t_1} \phi_1 *_{t_2} \dots *_{t_n} \phi_n$ such a path. Then
a strictly increasing reparametrization is of the form $\phi'_0
*_{t'_1} \phi_1 *_{t'_2} \dots *_{t'_n} \phi'_n$ where the $\phi'_i$
maps are strictly increasing reparametrizations of the $\phi_i$ maps.
Then the cocone $(\overline{f_i}:A_i \rightarrow A)$ is the unique
$U_\Delta$-final lift. The last sentence is a consequence of
\cite[Proposition~21.15]{topologicalcat}.  \epf

Let $f:X\rightarrow Y$ be a map of multipointed $d$-spaces. For every
$(\alpha,\beta)\in X^0\p X^0$, the set map
$\P^{top}f:\P^{top}_{\alpha,\beta} X \rightarrow
\P^{top}_{f^0(\alpha),f^0(\beta)} Y$ is continuous if
$\P^{top}_{\alpha,\beta} X$ and $\P^{top}_{f^0(\alpha),f^0(\beta)} Y$
are equipped with the Kelleyfication of the relative topology coming
from the inclusions $\P^{top}_{\alpha,\beta} X \subset
\ttop_\Delta([0,1],|X|)$ and $\P^{top}_{f^0(\alpha),f^0(\beta)} Y
\subset \ttop_\Delta([0,1],|Y|)$.

The following proposition will be used in the paper:

\bp \label{accessible} The functor $\P^{top}:\mdtop_\Delta \rightarrow
\top_\Delta$ is finitely accessible. \ep

\bpf Since the forgetful functor $\omega_\Delta:\top_\Delta
\rightarrow \set$ is topological by Proposition~\ref{bientop}, it
creates colimits by \cite[Proposition~21.15]{topologicalcat}. So it
suffices to prove that the functor $\omega_\Delta \circ
\P^{top}:\mdtop_\Delta \rightarrow \set$ is finitely accessible. The
latter fact is due to the construction of colimits in $\mdtop_\Delta$
(see Proposition~\ref{limcolim}).  \epf

Since the functor $\P^{top}:\mdtop_\Delta \rightarrow \top_\Delta$
preserves limits, it is a right adjoint by
\cite[Theorem~1.66]{MR95j:18001}. The left adjoint $G:\top_\Delta
\rightarrow \mdtop_\Delta$ is explicitly constructed in
Section~\ref{new}.

\section{Globular complexes}
\label{new}

Let $Z$ be a topological space. By definition, the \textit{topological
  globe of $Z$}, which is denoted by $\glob^{top}(Z)$, is the
multipointed $d$-space
\[
\left(\frac{\{\widehat{0},\widehat{1}\}\sqcup
    (Z\p_\Delta[0,1])}{(z,0)=(z',0)=\widehat{0},(z,1)=(z',1)=\widehat{1}},\{\widehat{0},\widehat{1}\},
  \P^{top} \glob^{top}(Z)\right)\] where $\P^{top} \glob^{top}(Z)$ is
the closure by strictly increasing reparametrizations of the set of
continuous maps $\{t\mapsto (z,t),z\in Z\}$. In particular,
$\glob^{top}(\varnothing)$ is the multipointed $d$-space
$\{\widehat{0},\widehat{1}\} =
(\{\widehat{0},\widehat{1}\},\{\widehat{0},\widehat{1}\},\varnothing)$.
Let \[\boxed{I^{gl,top} =
  \{\glob^{top}(\mathbf{S}^{n-1})\longrightarrow
  \glob^{top}(\mathbf{D}^{n}),n\geq 0\}}\] where the maps are the
image by $\glob(-)$ of the inclusions of spaces $\mathbf{S}^{n-1}
\subset \mathbf{D}^{n}$ with $n\geq 0$.

\bd A {\rm globular complex} is a multipointed $d$-space $X$ such that
the map $X^0 \rightarrow X$ is a relative $I^{gl,top}$-cell complex.
The set of cells of $X^0 \rightarrow X$ is called the {\rm globular
  decomposition} of the globular complex.  \ed

\begin{nota} $\vI^{top} = \glob^{top}(\mathbf{D}^0)$. \end{nota}

\bp \label{T1inc} Let $f:X\rightarrow Y$ be an element of
$\cell(I^{gl,top})$. Then $|f|:|X| \rightarrow |Y|$ is a Serre
cofibration of spaces. In particular, it is a closed $T_1$-inclusion
of spaces.  \ep

\bpf The continuous map $|\glob^{top}(\mathbf{S}^{n-1})| \rightarrow
|\glob^{top}(\mathbf{D}^{n})|$ is a cofibration by
\cite[Theorem~8.2]{4eme} (the fact that \cite{4eme} is written in the
category of compactly generated spaces does not matter here). Since
the class of cofibrations of the Quillen model structure of
$\top_\Delta$ is closed under pushout and transfinite composition, one
obtains the first sentence using Proposition~\ref{limcolim}. The map
$|f|:|X| \rightarrow |Y|$ is also a cofibration in $\mathcal{T}$ and
therefore a closed $T_1$-inclusion of spaces by
\cite[Lemma~2.4.5]{MR99h:55031}.  \epf

\bp \label{bonspace} Let $X$ be a globular complex.  Then the
topological space $|X|$ is compactly generated. \ep

\bpf By Proposition~\ref{limcolim} and Proposition~\ref{T1inc}, the
space $|X|$ is cofibrant.  Therefore it is weak Hausdorff. \epf

\bp \label{bonspace1} Let $X$ be a globular complex. Then the space
$\P^{top}X$ is compactly generated.  \ep

\bpf The space $\mathcal{T}_t([0,1],|X^{cof}|)$ is weak Hausdorff by
\cite[Lemma~5.2]{Ref_wH} since $|X^{cof}|$ is weak Hausdorff by
Proposition~\ref{bonspace}.  Thus the space
$\ttop_\Delta([0,1],|X^{cof}|)$ is weak Hausdorff since the
Kelleyfication functor adds open and closed subsets. So $\P^{top}X$
which is equipped with the Kelleyfication of the relative topology is
weak Hausdorff as well.  Thus it is compactly generated by
Proposition~\ref{inclusiontop}.  \epf

\bp \label{compact} (\cite[Lemma~1.1]{Ref_wH}) Let $X$ be a weak
Hausdorff space. Let $g:K\rightarrow X$ be a continuous map from a
compact $K$ to $X$. Then $g(K)$ is compact. \ep

\bp Let $X$ be a globular complex of $\mdtop_\Delta$. Then $X^0$ is a
discrete subspace of $X$. Let $\phi\in \P^{top}X$.  Then there exist
$0 = t_0 < \dots < t_n = 1$ such that
\begin{itemize}
\item $\phi(t_i)\in X^0$ for all $i\in \{0,1,\dots,n\}$ 
\item $\phi(]t_i,t_{i+1}[)\cap X^0=\varnothing$ for all $i\in \{0,1,\dots,n-1\}$ 
\item the restriction $\phi\res_{]t_i,t_{i+1}[}$ is
  one-to-one. 
\end{itemize}
\ep

\bpf By definition, the map $X^0 \rightarrow X$ is a relative
$I^{gl,top}$-cell complex. By Proposition~\ref{T1inc}, the continuous
map $X^0 \rightarrow |X|$ is a closed $T_1$-inclusion of spaces. So
$X^0$ is a discrete subspace of $X$. Since $[0,1]$ is compact and
since $|X|$ is compactly generated by Proposition~\ref{bonspace}, the
subset $\phi([0,1])$ is a closed subset and a compact subspace of
$|X|$ by Proposition~\ref{compact}.  Thus $\phi([0,1])\cap X^0$ is
finite. The rest of the statement is then clear.  \epf

Thus this new definition of globular complex coincides with the (very
long !) definition of \cite{model2}.

\begin{nota} The full subcategory of globular complexes of
  $\mdtop_\Delta$ is denoted by $\gltop$. \end{nota}

As an illustration of the objects of this section, it is now
constructed the left adjoint $G:\top_\Delta \rightarrow \mdtop_\Delta$
of the path space functor $\P^{top}:\mdtop_\Delta \rightarrow
\top_\Delta$.  Let $Z$ be a $\Delta$-generated space. If $Z$ is
non-empty connected, one has the natural set bijection
\[\mdtop_\Delta(\glob^{top}(Z),X) \iso \top_\Delta(Z,\P^{top}X)\] 
because of the cartesian closedness of $\top_\Delta$.  In the general
case, the $\Delta$-generated space $Z$ is homeomorphic to the disjoint
sum of its non-empty connected components by
Proposition~\ref{connexe}. Denote this situation by $Z\iso
\bigsqcup_{Z_i\in \pi_0(Z)} Z_i$. Then the functor $G$ is defined on
objects by
\[\boxed{G(Z)=\bigsqcup_{Z_i\in \pi_0(Z)}\glob^{top}(Z_i)}\] and in on obvious
way on morphisms.

\bp The functor $G:\top_\Delta \rightarrow \mdtop_\Delta$ is left
adjoint to the functor $\P^{top}:\mdtop_\Delta \rightarrow
\top_\Delta$. \ep 

\bpf It suffices to use the fact that $\top_\Delta$ is cartesian
closed.  \epf

\section{S-homotopy and strong deformation retract}
\label{retract}

The space $\mttop_\Delta((|X|,X^0),(|Y|,Y^0))$ is defined as the set
$\mtop_\Delta((|X|,X^0),(|Y|,Y^0))$ equipped with the Kelleyfication
of the relative topology coming from the set inclusion
\[\mtop_\Delta((|X|,X^0),(|Y|,Y^0)) \subset \ttop_\Delta(|X|,|Y|)\p \set(X^0,Y^0)\]
where $\set(X^0,Y^0)$ is equipped with the discrete topology. The
space \[\mdttop_\Delta((|X|,X^0,\P^{top}X),(|Y|,Y^0,\P^{top}Y))\] is
defined as the set
$\mdtop_\Delta((|X|,X^0,\P^{top}X),(|Y|,Y^0,\P^{top}Y))$ equipped with
the Kelleyfication of the relative topology coming from the set
inclusion
\[\mdtop_\Delta((|X|,X^0,\P^{top}X),(|Y|,Y^0,\P^{top}Y)) \subset
\mttop_\Delta((|X|,X^0),(|Y|,Y^0)).\]

\bd \label{S} (Compare with \cite[Definition~7.2 and
Proposition~7.5]{model3}) Let $f,g:X \rightrightarrows Y$ be two
morphisms of multipointed $d$-spaces. Then $f$ and $g$ are {\rm
  S-homotopic} if there exists a continuous map $H:[0,1] \rightarrow
\mdttop_\Delta(X,Y)$ called a {\rm S-homotopy} such that $H(0)=f$ and
$H(1)=g$. This situation is denoted by $\sim_S$. Two multipointed
$d$-spaces $X$ and $Y$ are {\rm S-homotopy equivalent} if and only if
there exist two morphisms $f:X\leftrightarrows Y:g$ with $g\circ
f\sim_S \id_X$ and $f\circ g\sim_S \id_Y$. \ed

\bd (Compare with \cite[Proposition~8.1]{model3}) Let $X$ be a
multipointed $d$-space. Let $U$ be a non-empty connected
$\Delta$-generated space. Then the multipointed $d$-space $U\boxtimes
X$ is defined as follows:
\begin{itemize}
\item Let $(U\boxtimes X)^0=X^0$.
\item Let $|U\boxtimes X|$ be defined by the pushout diagram of
  spaces
\[
\xymatrix{
  U \p_\Delta X^0 \fd{}\fr{} & U \p_\Delta |X| \fd{} \\
  X^0 \fr{} & \cocartesien {|U\boxtimes X|}.}
\]
\item The set $\P_{\alpha,\beta}^{top}(U\boxtimes X)$ is the smallest
  subset of continuous maps from $[0,1]$ to $|U\boxtimes X|$
  containing the continuous maps $\phi:[0,1] \rightarrow |U\boxtimes
  X|$ of the form $t\mapsto (u,\phi_2(t))$ where $u\in [0,1]$ and
  $\phi_2\in \P_{\alpha,\beta}^{top}X$ and closed under composition
  $*_{1/2}$ and strictly increasing reparametrization.
\end{itemize}
\ed 

\bd Let $X$ be a multipointed $d$-space. The multipointed $d$-space
$[0,1]\boxtimes X$ is called the {\rm cylinder} of $X$.  \ed

\bth \label{cyl} (Compare with \cite[Theorem~7.9]{model3}) Let $U$ be
a non-empty connected space. Let $X$ and $Y$ be two multipointed
$d$-spaces. Then there is a natural bijection of sets
\[\mdtop_\Delta(U\boxtimes X,Y) \iso
\top_\Delta(U,\mdttop_\Delta(X,Y)).\] \eth

\bpf Let $f:U\boxtimes X \rightarrow Y$ be a map of multipointed
$d$-spaces. Then by definition, one has a set map $f^0:X^0 \rightarrow
Y^0$ and a continuous map $|f|:|U\boxtimes X| \rightarrow |Y|$. By
adjunction, the composite $U\p_\Delta |X| \rightarrow |U\boxtimes X|
\rightarrow |Y|$ gives rise to a continuous map $\widehat{f}:U
\rightarrow \ttop_\Delta(|X|,|Y|)$. Let $\phi_2\in
\P_{\alpha,\beta}^{top}X$.  Then one has the equality
\[(\widehat{f}(u)\circ \phi_2)(t) = f(u,\phi_2(t))\] for any $(u,t)\in
U\p_\Delta [0,1]$. The continuous map $t\mapsto (u,\phi_2(t))$ is an
element of $\P_{\alpha,\beta}^{top}(U\boxtimes X)$ by definition of
the latter space. So the continuous map $t\mapsto f(u,\phi_2(t))$ is
an element of $\P_{f^0(\alpha),f^0(\beta)}^{top}Y$ since $f:U\boxtimes
X \rightarrow Y$ is a map of multipointed $d$-spaces. Thus for all
$u\in U$, $\widehat{f}(u)$ together with $f^0$ induce a map of
multipointed $d$-spaces from $X$ to $Y$. One obtains a set map from
$U$ to $\mdttop_\Delta(X,Y)$. It is continuous if and only if the
composite $U \longrightarrow \mdttop_\Delta(X,Y) \subset
\ttop_\Delta(|X|,|Y|)\p_\Delta \set(X^0,Y^0)$ is continuous by
definition of the topology of $\mdttop_\Delta(X,Y)$. By adjunction,
the latter map corresponds to the continuous map $(|f|,f^0)\in
\ttop_\Delta(U\p_\Delta |X|,|Y|)\p_\Delta \set(X^0,Y^0)$. Hence the
continuity and the commutative diagram of sets
\begin{equation}\label{1}
  \xymatrix{
    \mdtop_\Delta(U\boxtimes X,Y) \fr{f\mapsto (\widehat{f},f^0)}\fd{\subset}& \top_\Delta(U,\mdttop_\Delta(X,Y))\fd{\subset} \\
    \top_\Delta(U\p_\Delta |X|,|Y|)\p_\Delta \set(X^0,Y^0) \fr{\iso} & \top_\Delta(U,\ttop_\Delta(|X|,|Y|))\p_\Delta
    \set(X^0,Y^0)}
\end{equation}

Conversely, let $g:U \rightarrow \mdttop_\Delta(X,Y)$ be a continuous
map. The inclusion \[\mdttop_\Delta(X,Y) \subset
\ttop_\Delta(|X|,|Y|)\p_\Delta \set(X^0,Y^0)\] gives rise to a set map
$\widetilde{g}^0:X^0 \rightarrow Y^0$ since $U$ is connected non-empty
and since $\set(X^0,Y^0)$ is discrete and to a continuous map
$|\widetilde{g}|:|U\boxtimes X| \rightarrow |Y|$ by adjunction. Let
$\phi\in \P_{\alpha,\beta}^{top}(U\boxtimes X)$.  By construction,
$\phi$ is a composition of continuous paths of the form
$(u,\phi_2(-))$ with $\phi_2\in \P_{\alpha,\beta}^{top}X$ and $u\in
U$.  Then $|\widetilde{g}|(u,\phi_2(t))=g(u)\circ \phi_2(t)$.  Since
$g(u)\in \mdtop_\Delta(X,Y)$, $|\widetilde{g}|\circ \phi
\in\P_{\widetilde{g}^0(\alpha),\widetilde{g}^0(\beta)}^{top}Y$ and
therefore $\widetilde{g}\in \mdtop_\Delta(U\boxtimes X,Y)$.  Hence the
commutative diagram of sets
\begin{equation}\label{2}
  \xymatrix{
    \top_\Delta(U,\mdttop_\Delta(X,Y))\fd{\subset}\fr{g\mapsto\widetilde{g}} &\mdtop_\Delta(U\boxtimes X,Y) \fd{\subset} \\
    \top_\Delta(U,\ttop_\Delta(|X|,|Y|))\p_\Delta \set(X^0,Y^0)\fr{\iso} &\top_\Delta(U\p_\Delta |X|,|Y|)\p_\Delta \set(X^0,Y^0)  }
\end{equation}

The commutativity of the diagrams $(\ref{1})$ and $(\ref{2})$ implies
that the set maps $f\mapsto (\widehat{f},f^0)$ and
$g\mapsto\widetilde{g}$ are inverse to each other, hence the result.
\epf

\bp \label{passage} Let $f,g:X \rightrightarrows Y$ be two S-homotopic
morphisms of multipointed $d$-spaces. Then $f^0=g^0$ and the
continuous maps $|f|,|g|:|X| \rightrightarrows |Y|$ are homotopic.
Moreover, for all $(\alpha,\beta)\in X^0\p X^0$, the pair of
continuous maps $\P^{top}f,\P^{top}g:\P^{top}_{\alpha,\beta} X
\rightrightarrows \P^{top}_{f^0(\alpha),f^0(\beta)} Y$ are homotopic.
\ep

\bpf Let $H:[0,1] \rightarrow \mdttop_\Delta(X,Y)$ be a S-homotopy between
$f$ and $g$. Then the composite 
\[[0,1]\p_\Delta |X| \longrightarrow |[0,1]\boxtimes X|
\stackrel{\widetilde{H}}\longrightarrow |Y|\] gives a homotopy between
$|f|$ and $|g|$. The mapping $t\mapsto \P^{top}(H(t))$ induces a set
map from $[0,1]$ to $\ttop_\Delta(\P^{top}_{\alpha,\beta}
X,\P^{top}_{f^0(\alpha),f^0(\beta)} Y)$ for any $(\alpha,\beta)\in
X^0\p X^0$. The latter set map is continuous since it corresponds by
adjunction to the continuous mapping $(t,\phi)\mapsto |H(t)|\circ
\phi$ from $[0,1]\p_\Delta \P^{top}_{\alpha,\beta} X$ to
$\P^{top}_{f^0(\alpha),f^0(\beta)} Y$.  So one obtains a homotopy
between $\P^{top} f$ and $\P^{top}g$.  Finally, the composite \[[0,1]
\stackrel{H}\longrightarrow \mdttop_\Delta(X,Y) \subset
\ttop_\Delta(|X|,|Y|) \p_\Delta \set(X^0,Y^0) \longrightarrow
\set(X^0,Y^0)\] is constant since $[0,1]$ is connected. So $f^0=g^0$.
\epf

\bd A map $i:A\rightarrow B$ of multipointed $d$-spaces is an {\rm
  inclusion of a strong deformation retract} if there exists a
continuous map $H:[0,1] \rightarrow \mdttop_\Delta(B,B)$ called the
{\rm deformation retract} such that
\begin{itemize}
\item $H(0)=\id_B$, $H(1)=i\circ r$ where $r:B\rightarrow A$ is a map
  of multipointed $d$-spaces such that $r\circ i=\id_A$
\item $|\widetilde{H}|(t,i(a))=i(a)$ for all $a\in A$ and for all
  $t\in [0,1]$.
\end{itemize}
\ed 

\bd A continuous map $i:A\rightarrow B$ of general topological spaces
is an {\rm inclusion of a strong deformation retract of spaces} if
there exists a continuous map $H:[0,1]\p B \rightarrow B$ called the
{\rm deformation retract} such that
\begin{itemize}
\item $H(0,-)=\id_B$, $H(1,-)=i\circ r$ where $r:B\rightarrow A$ is a map
  of multipointed $d$-spaces such that $r\circ i=\id_A$
\item $H(t,i(a))=i(a)$ for all $a\in A$ and for all $t\in [0,1]$.
\end{itemize}
\ed 

\bp \label{passage2} Let $i:A \rightarrow B$ be an inclusion of a
strong deformation retract of multipointed $d$-spaces. Then
\begin{itemize}
\item The set map $i^0:A^0\rightarrow B^0$ is bijective.
\item The continuous map $|i|:|A| \rightarrow |B|$ is an inclusion of
  a strong deformation retract of spaces.
\item The continuous map $\P^{top}i:\P^{top}_{\alpha,\beta} A
  \rightarrow \P^{top}_{f^0(\alpha),f^0(\beta)} B$ is an inclusion of
  a strong deformation retract of spaces for all $(\alpha,\beta)\in
  A^0\p A^0$.
\end{itemize} 
\ep

\bpf Consequence of Proposition~\ref{passage}.  \epf

\bth \label{push} The class of inclusions of a strong deformation
retract of multipointed $d$-spaces is closed under pushout. \eth

\bpf We mimic the proof given in
\cite[Proposition~2.4.9]{MR99h:55031}.  Consider a pushout diagram of
multipointed $d$-spaces
\[
\xymatrix{
A \fr{}\fd{i} & C \fd{j} \\
B \fr{} & \cocartesien D}
\]
where $i$ is an inclusion of a strong deformation retract. Let
$K:[0,1] \rightarrow \mdttop_\Delta(B,B)$ be the corresponding
deformation retract and let $r:B\rightarrow A$ be the retraction. The
commutative diagram of multipointed $d$-spaces
\[
\xymatrix{
A \ar@{=}[rrrd]\fr{}\fd{i} & C \ar@{->}[d]{} \ar@{=}[rrrd] & &\\
B \ar@{->}[r] \ar@{->}[rrrd]_-{r}& \cocartesien D \ar@{-->}[rrrd]&& A \fr{}\ar@{=}[d] & C \fd{}\\
&&& A \fr{} & \cocartesien C}
\]
gives the retraction $s:D\rightarrow C$. Denote by $\widetilde{\id_C}$
the map corresponding to the constant map from $[0,1]$ to
$\mdttop_\Delta(C,C)$ taking any element of $[0,1]$ to the identity of
$C$.  Then consider the commutative diagram of multipointed $d$-spaces
\[
\xymatrix{
[0,1]\boxtimes A \fr{}\fd{} & [0,1]\boxtimes C \ar@{->}[d]{} \ar@{->}[rrdd]^-{\widetilde{\id_C}} & &\\
[0,1]\boxtimes B \ar@{->}[r] \ar@{->}[rrdd]_-{\widetilde{K}}& \cocartesien [0,1]\boxtimes D \ar@{-->}[rrdd]&& \\
&& & C \fd{}\\
&& B \fr{} & D}
\]
Since the functor $[0,1]\boxtimes -$ preserves colimits by
Theorem~\ref{cyl}, the universal property of the pushout induces a map
of multipointed $d$-spaces $\widetilde{H} : [0,1]\boxtimes D
\rightarrow D$. Denote by $H$ the corresponding continuous map from
$[0,1]$ to $\mdttop_\Delta(D,D)$. By Proposition~\ref{limcolim}, one
obtains the commutative diagram of topological spaces
\[
\xymatrix{
  |[0,1]\boxtimes A| \fr{}\fd{} & |[0,1]\boxtimes C| \ar@{->}[d]{} \ar@{->}[rrdd]^-{|\widetilde{\id_C}|} & &\\
  |[0,1]\boxtimes B| \ar@{->}[r] \ar@{->}[rrdd]_-{|\widetilde{K}|}& \cocartesien {|[0,1]\boxtimes D|} \ar@{-->}[rrdd]&& \\
  && & |C| \fd{}\\
  && |B| \fr{} & |D|.}
\]
By construction, one then has $|\widetilde{H}|(t,c) = j(c)$ for all
$c\in C$ and all $t\in[0,1]$ and $|\widetilde{H}|(0,d) = d$ for all
$d\in D$. Since $|\widetilde{K}|(1,b)\in |i|(|A|)$ for all $b\in B$,
it follows that $|\widetilde{H}|(1,d)\in |j|(|C|)$ for all $d\in D$.
Since $|j|$ is an inclusion of spaces by Proposition~\ref{passage2},
$H$ is a deformation retract as required.  \epf

\section{The combinatorial model structure}
\label{comb}

\bd (Compare with \cite[Definition~11.6]{model3}) A morphism
$f:X\rightarrow Y$ of multipointed $d$-spaces is a {\rm weak
  S-homotopy equivalence} if $f^0:X^0 \rightarrow Y^0$ is a bijection
of sets and if $\P^{top} f:\P^{top}_{\alpha,\beta} X\rightarrow
\P^{top}_{f^0(\alpha),f^0(\beta)} Y$ is a weak homotopy equivalence of
topological spaces.  The class of weak S-homotopy equivalences of
multipointed $d$-spaces is denoted by $\mathcal{W}^{top}$. \ed

If $C : \varnothing \rightarrow \{0\}$ and $R : \{0,1\} \rightarrow
\{0\}$ are set maps viewed as morphisms of multipointed $d$-spaces,
let $I^{gl,top}_{+} = I^{gl,top} \cup \{C,R\}$. Let
\[\boxed{J^{gl,top} =
  \{\glob^{top}(\mathbf{D}^{n}\p_\Delta\{0\})\longrightarrow
  \glob^{top}(\mathbf{D}^{n}\p_\Delta[0,1]),n\geq 0\}}\]
where the maps are the image by $\glob(-)$ of the inclusions of spaces
$\mathbf{D}^{n}\p_\Delta\{0\} \subset \mathbf{D}^{n}\p_\Delta[0,1]$
with $n\geq 0$.

\bp \label{prediag} Let $f:X \rightarrow Y$ be a morphism of
multipointed $d$-spaces. Then $f$ satisfies the right lifting property
with respect to $\{C,R\}$ if and only if $f^0:X^0 \rightarrow Y^0$ is
bijective. \ep

\bpf See \cite[Proposition~16.2]{model3} or \cite[Lemma~4.4
(5)]{nonexistence}.  \epf

\bp \label{diag} (Compare with \cite[Proposition~13.2]{model3}) Let
$f:X\rightarrow Y$ be a morphism of multipointed $d$-spaces.  Let
$g:U\rightarrow V$ be a continuous map.  Then $f$ satisfies the right
lifting property with respect to $\glob^{top}(g)$ if and only if for
all $(\alpha,\beta)\in X^0\p X^0$, $\P^{top}f:
\P^{top}_{\alpha,\beta}X \rightarrow
\P^{top}_{f^0(\alpha),f^0(\beta)}Y$ satisfies the right lifting
property with respect to $g$.  \ep

\bpf Consider a commutative diagram of solid arrows of spaces
\[
\xymatrix{
U \fd{g}\fr{} & \P^{top}_{\alpha,\beta}X \ar@{->}[d]^-{\P^{top}f} \\
V \fr{} \ar@{-->}[ru]^-{\ell}& \P^{top}_{f^0(\alpha),f^0(\beta)}Y.}
\]
Since $\top_\Delta$ is cartesian closed by
Proposition~\ref{cartesianclosed}, the existence of the lift $\ell$ is
equivalent to the existence of the lift $\ell'$ in the commutative
diagram of solid arrows of multipointed spaces
\[
\xymatrix{
\glob^{top}(U) \fd{g}\fr{} & X \fd{f} \\
\glob^{top}(V) \fr{} \ar@{-->}[ru]^-{\ell'}& Y.}
\]
\epf

\bp \label{passage3} Let $i:A\rightarrow B$ be an inclusion of a
strong deformation retract of spaces between $\Delta$-generated
spaces. Then $\glob^{top}(i):\glob^{top}(A) \rightarrow
\glob^{top}(B)$ is an inclusion of a strong deformation retract of
$\mdtop_\Delta$.  \ep

\bpf Let $r:B\rightarrow A$ be the corresponding retraction and let
$H:[0,1]\p_\Delta B\rightarrow B$ be the corresponding deformation
retract.  Then $\glob^{top}(r):\glob^{top}(B) \rightarrow
\glob^{top}(A)$ is the corresponding retraction of multipointed
$d$-spaces and $\glob^{top}(H):\glob^{top}([0,1]\p_\Delta B)
\rightarrow \glob^{top}(B)$ is the corresponding deformation retract
of multipointed $d$-spaces since there is an isomorphism of
multipointed $d$-spaces $\glob^{top}([0,1]\p_\Delta B) \iso
[0,1]\boxtimes \glob^{top}(B)$.  \epf

\bth \label{model} There exists a unique cofibrantly generated model
structure on $\mdtop_\Delta$ such that the weak equivalences are the
weak S-homotopy equivalences, such that $I^{gl,top}_{+}=I^{gl,top}\cup
\{C,R\}$ is the set of generating cofibrations and such that
$J^{gl,top}$ is the set of trivial generating cofibrations.  Moreover,
one has:
\begin{enumerate}
\item The cellular objects are exactly the globular complexes.
\item A map $f:X\rightarrow Y$ of multipointed $d$-spaces is a (resp.
  trivial) fibration if and only if for all $(\alpha,\beta)\in X^0\p
  X^0$, the continuous map $\P^{top}_{\alpha,\beta}X \rightarrow
  \P^{top}_{f^0(\alpha),f^0(\beta)}Y$ is a (resp. trivial) fibration
  of spaces.
\item Every multipointed $d$-space is fibrant.
\item The model structure is right proper and simplicial.
\item Two cofibrant multipointed $d$-spaces are weakly S-homotopy
  equivalent if and only they are S-homotopy equivalent.
\end{enumerate}
\eth

\bpf We have to prove the following facts
\cite[Theorem~2.1.19]{MR99h:55031}:
\begin{enumerate}
\item The class $\mathcal{W}^{top}$ satisfies the two-out-of-three property
  and is closed under retracts.
\item The domains of $I^{gl,top}_{+}$ are small relative to
  $\cell(I^{gl,top}_{+})$.
\item The domains of $J^{gl,top}$ are small relative to
  $\cell(J^{gl,top})$.
\item $\cell(J^{gl,top}) \subset \mathcal{W}^{top} \cap
  \cof(I^{gl,top}_{+})$. 
\item $\inj(I^{gl,top}_{+})=\mathcal{W}^{top} \cap \inj(J^{gl,top})$.
\end{enumerate}
That $\mathcal{W}^{top}$ satisfies the two-out-of-three property and
is closed under retracts is clear. By
\cite[Proposition~1.3]{MR1780498}, every object of $\mdtop_\Delta$ is
small relative to the whole class of morphisms since $\mdtop_\Delta$
is locally presentable by Theorem~\ref{locpre}. Hence Assertions (2)
and (3). By Proposition~\ref{prediag} and Proposition~\ref{diag}, a
map of multipointed $d$-spaces $f:X \rightarrow Y$ belongs to
$\inj(I^{gl,top}_{+})$ if and only if the set map $f^0:X^0 \rightarrow
Y^0$ is bijective and for all $(\alpha,\beta)\in X^0\p X^0$, the map
$\P^{top}_{\alpha,\beta}X \rightarrow \P^{top}_{f(\alpha),f(\beta)}Y$
is a trivial fibration of spaces. By Proposition~\ref{diag} again, a
map of multipointed $d$-spaces $f:X\rightarrow Y$ belongs to
$\inj(J^{gl,top})$ if and only if for all $(\alpha,\beta)\in X^0\p
X^0$, the map $\P^{top}_{\alpha,\beta}X \rightarrow
\P^{top}_{f(\alpha),f(\beta)}Y$ is a fibration of spaces. Hence the
fifth assertion.

It remains to prove the fourth assertion which is the most delicate
part of the proof. The set inclusion $\cell(J^{gl,top}) \subset
\cell(I^{gl,top})$ comes from the obvious set inclusion $J^{gl,top}
\subset \cof(I^{gl,top})$.  It remains to prove the set inclusion
$\cell(J^{gl,top}) \subset \mathcal{W}^{top}$.  First of all, consider
the case of a map $f:X\rightarrow Y$ which is the pushout of one
element of $J^{gl,top}$ (for some $n\geq 0$):
\[
\xymatrix{
\glob^{top}(\mathbf{D}^n\p_\Delta \{0\}) \fd{}\fr{} & X \fd{} \\
\glob^{top}(\mathbf{D}^n\p_\Delta [0,1]) \fr{} & \cocartesien Y.}
\]
Since each map of $J^{gl,top}$ is an inclusion of a strong deformation
retract of multipointed $d$-spaces by Proposition~\ref{passage3}, the
map $f:X\rightarrow Y$ is itself an inclusion of a strong deformation
retract by Theorem~\ref{push}.  So by Proposition~\ref{passage2}, the
continuous map $\P^{top}f : \P^{top}X \rightarrow \P^{top}Y$ is an
inclusion of a strong deformation retract of spaces.  By
Proposition~\ref{T1inc}, the subset $|Y|\backslash |X|$ is open in
$|Y|$. Let $\phi\in \P^{top}Y \backslash \P^{top}X$. By
Proposition~\ref{limcolim}, one has $\phi=\phi_0 *_{t_1} \phi_1
*_{t_2} \dots *_{t_n} \phi_n$ with the $\phi_i$ being execution paths
of the images of $\glob^{top}(\mathbf{D}^n\p_\Delta [0,1])$ or of $X$
in $Y$.  Since $\phi\notin \P^{top}X$, the set
$\phi^{-1}(|Y|\backslash |X|)$ is a non-empty open subset of $[0,1]$.
Let $[a,b] \subset \phi^{-1}(|Y|\backslash |X|)$ with $0 \leq a < b
\leq 1$. Then $\phi$ belongs to $N([a,b],|Y|\backslash |X|)$ which is
an open of the compact-open topology. Since the Kelleyfication functor
adds open subsets, the set $N([a,b],|Y|\backslash |X|)$ is an open
neighborhood of $\phi$ included in $\P^{top}Y \backslash \P^{top}X$.
So $\P^{top}Y \backslash \P^{top}X$ is open in $\P^{top}Y$. Therefore,
the continuous map $\P^{top}f : \P^{top}X \rightarrow \P^{top}Y$ is a
closed inclusion of a strong deformation retract of
spaces~\footnote{This step of the proof is necessary. See
  Proposition~\ref{ir} of this paper.}.  This class of maps is closed
under transfinite composition in $\mathcal{T}$, and in $\top_\Delta$,
since it is the class of trivial cofibrations of the Str{\o}m model
structure on the category of general topological spaces
\cite{MR35:2284} \cite{MR39:4846} \cite{ruse}. So for any map $f\in
\cell(J^{gl,top})$, the continuous map $\P^{top}f : \P^{top}X
\rightarrow \P^{top}Y$ is a closed inclusion of a strong deformation
retract of spaces since the functor $\P^{top}$ is finitely accessible
by Proposition~\ref{accessible}.  Therefore it is a weak homotopy
equivalence.  Hence the set inclusion $\cell(J^{gl,top}) \subset
\mathcal{W}^{top}$.

It remains to prove the last five assertions of the statement of the
theorem. The first one is obvious. The second one is explained above.
The third one comes from the second one and from the fact that all
topological spaces are fibrant (the terminal object of $\mdtop_\Delta$
is $(\{0\},\{0\},\{0\})$). The model structure is right proper since
all objects are fibrant. The construction of the simplicial structure
of this model category is postponed until Appendix~\ref{simp}. The
very last assertion is due to the construction of the simplicial
structure and to the fact that $X\ot \Delta[1]\iso [0,1]\boxtimes X$
since $[0,1]$ is non-empty and connected. \epf

Here are some comments about cofibrancy. The set
$\{0\}=(\{0\},\{0\},\varnothing)$ is cofibrant. More generally, any
set $E=(E,E,\varnothing)$ is cofibrant. The multipointed $d$-space
$(\{0\},\varnothing,\varnothing)$ is not cofibrant. Indeed, its
cofibrant replacement is equal to the initial multipointed $d$-space
$\varnothing=(\varnothing, \varnothing, \varnothing)$.  Finally, the
terminal multipointed $d$-space $(\{0\},\{0\},\{0\})$ is not
cofibrant.

\begin{conj} The model category $\mdtop_\Delta$ is left proper.
\end{conj}

\section{Comparing multipointed $d$-spaces and flows}
\label{comp}

We describe a functor $cat$ from the category of multipointed
$d$-spaces to the category of flows which generalizes the functor
constructed in \cite{model2}.

\bd \cite{model3} A {\rm flow} $X$ is a small category without
identity maps enriched over topological spaces. The composition law of
a flow is denoted by $*$.  The set of objects is denoted by $X^0$.
The space of morphisms from $\alpha$ to $\beta$ is denoted by
$\P_{\alpha,\beta} X$. Let $\P X$ be the disjoint sum of the spaces
$\P_{\alpha,\beta} X$. A morphism of flows $f:X \rightarrow Y$ is a
set map $f^0:X^0 \rightarrow Y^0$ together with a continuous map $\P
f:\P X \rightarrow \P Y$ preserving the structure.  The corresponding
category is denoted by $\dtop(\top_\Delta)$,
$\dtop(\top_\mathcal{K})$, etc... depending on the category of
topological spaces we are considering. \ed

Let $Z$ be a topological space. The flow $\glob(Z)$ is defined by
\begin{itemize}
\item $\glob(Z)^0=\{\widehat{0},\widehat{1}\}$, 
\item $\P \glob(Z)= \P_{\widehat{0},\widehat{1}} \glob(Z) = Z$,
\item $s=\widehat{0}$, $t=\widehat{1}$ and a trivial composition law.
\end{itemize}
It is called the \textit{globe} of the space $Z$.

Let $X$ be a multipointed $d$-space. Consider for every
$(\alpha,\beta)\in X^0 \p X^0$ the coequalizer of sets
\[\P_{\alpha,\beta}X = \liminj\left( \P^{top}_{\alpha,\beta}X\p_\Delta
  \P^{top}\vI^{top} \rightrightarrows
  \P^{top}_{\alpha,\beta}X\right)\] where the two maps are
$(c,\phi)\mapsto c \circ \id_{\vI^{top}} = c$ and $(c,\phi)\mapsto c
\circ \phi$. Let $[-]_{\alpha,\beta}:\P^{top}_{\alpha,\beta}X
\rightarrow \P_{\alpha,\beta}X$ be the canonical set map.  The set
$\P_{\alpha,\beta}X$ is equipped with the final topology.

\bth Let $X$ be a multipointed $d$-space. Then there exists a flow
$cat_\Delta(X)$ with $cat_\Delta(X)^0=X^0$,
$\P_{\alpha,\beta}cat_\Delta(X)= \P_{\alpha,\beta}X$ and the
composition law $*:\P_{\alpha,\beta}X \p_\Delta \P_{\beta,\gamma}X
\rightarrow \P_{\alpha,\gamma}X$ is for every triple
$(\alpha,\beta,\gamma)\in X^0\p X^0\p X^0$ the unique map making the
following diagram commutative:
\[
\xymatrix{
\P_{\alpha,\beta}^{top}X \p_\Delta \P_{\beta,\gamma}^{top}X
\fr{}\fd{[-]_{\alpha,\beta}\p_\Delta [-]_{\beta,\gamma}} &
\P_{\alpha,\gamma}^{top}X \fd{[-]_{\alpha,\gamma}} \\
\P_{\alpha,\beta}X \p_\Delta \P_{\beta,\gamma}X \fr{} &
\P_{\alpha,\gamma}X}
\] 
where the map $\P_{\alpha,\beta}^{top}X \p_\Delta
\P_{\beta,\gamma}^{top}X \rightarrow \P_{\alpha,\gamma}^{top}X$ is the
continuous map defined by the concatenation of continuous paths:
$(c_1,c_2)\in \P_{\alpha,\beta}^{top}X \p_\Delta
\P_{\beta,\gamma}^{top}X$ is sent to $c_1*_{1/2}c_2$.  The mapping $X
\mapsto cat_\Delta(X)$ induces a functor from $\mdtop_\Delta$ to
$\dtop(\top_\Delta)$.  \eth

\bpf The existence and the uniqueness of the continuous map
$\P_{\alpha,\beta}X \p_\Delta \P_{\beta,\gamma}X \rightarrow
\P_{\alpha,\gamma}X$ for any triple $(\alpha,\beta,\gamma)\in X^0\p
X^0 \p X^0$ comes from the fact that the topological space
$\P_{\alpha,\beta}X$ is the quotient of the set
$\P^{top}_{\alpha,\beta}X$ by the equivalence relation generated by
the identifications $c \circ \phi = c$ equipped with the final
topology. This defines a strictly associative law because the
concatenation of continuous paths is associative up to strictly
increasing reparametrization. The functoriality is obvious. \epf

\bp The functor $cat_\Delta:\mdtop_\Delta \rightarrow
\dtop(\top_\Delta)$ does not have any right adjoint. \ep

\bpf Suppose that the right adjoint exists. Let $Z$ be a
$\Delta$-generated space. Then for any multipointed $d$-space $X$, one
would have the natural isomorphism \[\mdtop_\Delta(X,Y) \iso
\dtop(\top_\Delta)(cat_\Delta(X),\glob(Z))\] for some multipointed
$d$-space $Y$.

With $X=X^0$, one obtains the equalities \[\mdtop_\Delta(X,Y) =
\set(X^0,Y^0)\] and \[\dtop(\top_\Delta)(cat_\Delta(X),\glob(Z)) =
\set(X^0,\glob(Z)^0).\] By an application of Yoneda's lemma within the
category of sets, one deduces that $Y^0 =
\{\widehat{0},\widehat{1}\}$.

Consider now the case $X = \vI^{top}$. The set
$\P_{\alpha,\beta}^{top}Y$ is non-empty if and only if there exists a
map $f:\vI^{top} \rightarrow Y$ with $f^0(\widehat{0}) = \alpha$ and
$f^0(\widehat{1}) = \beta$. The map $f:\vI^{top} \rightarrow Y$
corresponds by adjunction to a map $\vI=\glob(\{0\}) \rightarrow
\glob(Z)$. Thus one obtains $\P^{top}Y =
\P_{\{\widehat{0},\widehat{1}\}}^{top}Y$.

Consider now the case $X = \glob^{top}(T)$ for some $\Delta$-generated
space $T$. Then one has the isomorphisms $\mdtop_\Delta(X,Y)\iso
\top_\Delta(T,\P_{\{\widehat{0},\widehat{1}\}}^{top}Y)$ and
$\dtop(\top_\Delta)(cat_\Delta(X),\glob(Z)) \iso \top_\Delta(T,Z)$. So
by Yoneda's lemma applied within the category $\top_\Delta$, one
obtains the isomorphism $\P^{top}Y =
\P_{\{\widehat{0},\widehat{1}\}}^{top}Y \iso Z$.

Take $Z=\{0\}$. Then $Y$ is a multipointed $d$-space with unique
initial state $\widehat{0}$, with unique final state $\widehat{1}$ and
with $\P^{top}Y = \P_{\{\widehat{0},\widehat{1}\}}^{top}Y =\{0\} \neq
\varnothing$. Thus there exists a continuous map $\phi:[0,1]
\rightarrow |Y|$ with $\phi(0)=\widehat{0}$, $\phi(1)=\widehat{1}$
which is an execution path of $Y$. So any strictly increasing
reparametrization of $\phi$ is an execution path of $Y$. So
$\widehat{0}=\widehat{1}$ in $Y^0$.  Contradiction. \epf

\begin{cor} The functor $cat_\Delta : \mdtop_\Delta \longrightarrow
  \dtop_\Delta$ is not colimit-preserving. \end{cor}

\bpf The category $\mdtop_\Delta$ is a topological fibre-small
category over the category of sets. Therefore, it satisfies the
hypothesis of the opposite of the special adjoint functor
theorem. Thus, the functor $cat_\Delta$ is colimit-preserving if and
only if it has a right adjoint. \epf

The colimit-preserving functor $H\circ i_\Delta^\mathcal{K} :
\top_\Delta \rightarrow \cgtop$ induces a colimit-preserving functor
$\dtop(\top_\Delta) \rightarrow \dtop(\cgtop)$ from the category of
flows enriched over $\Delta$-generated spaces to that of flows
enriched over compactly generated topological spaces (the latter
category is exactly the category used in \cite{model3} and in
\cite{model2}). So the composite functor
\[cat:\mdtop_\Delta \stackrel{cat_\Delta} \longrightarrow
\dtop(\top_\Delta) \longrightarrow \dtop(\cgtop)\] coincides on globular
complexes with the functor constructed in \cite{model2} since the
globular complexes of \cite{model2} are exactly those defined in
Section~\ref{new} and since one has the equality $cat(\glob^{top}(Z))
= \glob(Z)$ for any compactly generated space $Z$.

\bth \label{equivalence} The composite functor 
\[\mdtop_\Delta \stackrel{(-)^{cof}}\longrightarrow \gltop
\stackrel{cat}\longrightarrow \ho(\dtop(\cgtop))\] induces an
equivalence of categories $\ho(\mdtop_\Delta) \simeq \ho(\dtop(\cgtop))$
between the category of multipointed $d$-spaces up to weak S-homotopy
and the category of flows up to weak S-homotopy where a map of flows
$f:X\rightarrow Y$ is a weak S-homotopy if and only if $f^0:X^0
\rightarrow Y^0$ is a bijection and $\P f:\P X \rightarrow \P Y$ is a
weak homotopy equivalence. \eth

\bpf By \cite[Theorem~V.4.2]{model2}, the restriction of the functor
cat to the category $\gltop$ of globular complexes induces a
categorical equivalence $\gltop[(\mathcal{W}^{top})^{-1}] \simeq
\ho(\dtop(\cgtop))$ where $\gltop[(\mathcal{W}^{top})^{-1}]$ is the
categorical localization of $\gltop$ by the class of weak S-homotopy
equivalences $\mathcal{W}^{top}$.  Every multipointed $d$-space is
fibrant. Therefore the localization $\gltop[(\mathcal{W}^{top})^{-1}]$
is isomorphic to the quotient $\gltop/\!\sim_S$ where $\sim_S$ is the
congruence relation on morphisms induced by S-homotopy (see
Definition~\ref{S}).  The inclusion functor $\gltop\subset
(\mdtop_\Delta)_{cof}$ from the category of globular complexes to that
of cofibrant multipointed $d$-spaces induces a full and faithful
functor $\gltop/\!\sim_S \rightarrow (\mdtop_\Delta)_{cof}/\!\sim_S$.
The right-hand category is equivalent to $\ho(\mdtop_\Delta)$ by
\cite[Section~7.5.6]{ref_model2} or
\cite[Proposition~1.2.3]{MR99h:55031}.  \epf

Recall that there exists a unique model structure on $\dtop(\cgtop)$
such that \cite{model3}
\begin{itemize}
\item The set of generating cofibrations is $I^{gl}_+ =
  I^{gl}\cup\{C,R\}$ with \[\boxed{I^{gl} =
    \{\glob(\mathbf{S}^{n-1})\longrightarrow
    \glob(\mathbf{D}^{n}),n\geq 0\}}\] where the maps are the image by
  $\glob(-)$ of the inclusions of spaces $\mathbf{S}^{n-1} \subset
  \mathbf{D}^{n}$ with $n\geq 0$.
\item The set of generating trivial cofibration is \[\boxed{J^{gl} =
    \{\glob(\mathbf{D}^{n}\p\{0\})\longrightarrow
    \glob(\mathbf{D}^{n}\p[0,1]),n\geq 0\}}\] where the maps are the
  image by $\glob(-)$ of the inclusions of spaces
  $\mathbf{D}^{n}\p\{0\} \subset \mathbf{D}^{n}\p[0,1]$ with $n\geq
  0$.
\item The weak equivalences are the weak S-homotopy equivalences.
\end{itemize}

\bp \label{path-cofibrant} Let $X$ be a cofibrant flow of
$\dtop(\cgtop)$. Then the path space $\P X$ is cofibrant. \ep

\bpf[Sketch of proof] The cofibrant flow $X$ is a retract of a
cellular flow $\overline{X}$ and the space $\P X$ is then a retract of
the space $\P \overline{X}$. Thus one can suppose $X$ cellular, i.e.
$\varnothing \rightarrow X$ belonging to $\cell(I^{gl}_+)$. By \cite[
Proposition~7.1]{2eme} (see also Proposition~\ref{calcul_explicite} of
this paper and \cite[Proposition~15.1]{model3}), the continuous map
$\varnothing \rightarrow \P X$ is a transfinite composition of a
$\lambda$-sequence $Y:\lambda \rightarrow \cgtop$ such that for any
$\mu<\lambda$, the map $Y_\mu \rightarrow Y_{\mu+1}$ is a pushout of a
map of the form $\id_{X_1} \p \dots \p i_n \p \dots \p \id_{X_p}$ with
$p\geq 0$ and with $i_n:\mathbf{S}^{n-1} \subset \mathbf{D}^n$ being
the usual inclusion with $n\geq 0$. Suppose that the set of
ordinals \[\{\mu\leq \lambda, Y_\mu \hbox{ is not cofibrant}\}\] is
non-empty.  Consider the smallest element $\mu_0$. Then the continuous
map $\varnothing \rightarrow \P Y_{\mu_0}$ is a transfinite
composition of pushout of maps of the form $\id_{X_1} \p \dots \p i_n
\p \dots \p \id_{X_p}$ with all spaces $X_i$ cofibrant since the
spaces $X_i$ are built using the path spaces of the flows $Y_\mu$ with
$\mu<\mu_0$.  So the space $\P Y_{\mu_0}$ is cofibrant because the
model category $\top_\Delta$ is monoidal. Contradiction.  \epf

Note that the analogous statement for multipointed $d$-spaces is
false. Indeed, the multipointed $d$-space $\vI^{top}$ is cofibrant
since $\vI^{top} \iso \glob^{top}(\mathbf{D}^0)$.  And the space
$\P^{top} \vI^{top}$ is the space of strictly increasing continuous
maps from $[0,1]$ to itself preserving $0$ and $1$. The latter space
is not cofibrant.

Let us conclude the comparison of multipointed $d$-spaces and flows
by:

\bth \label{QA} There exists a unique model structure on
$\dtop(\top_\Delta)$ such that the set of generating (resp. trivial)
cofibrations is $I^{gl}_+$ (resp. $J^{gl}$) and such that the weak
equivalences are the weak S-homotopy equivalences. Moreover:
\begin{enumerate}
\item The model category $\dtop(\top_\Delta)$ is proper simplicial and
  combinatorial.
\item The categorical adjunction $\dtop(\top_\Delta) \leftrightarrows
  \dtop(\cgtop)$ is a Quillen equivalence.
\item The functor $cat_\Delta:\mdtop_\Delta \rightarrow
  \dtop(\top_\Delta)$ preserves cofibrations, trivial cofibrations and
  weak S-homotopy equivalences between cofibrant objects.
\end{enumerate}
\eth 

\bpf[Sketch of proof] The construction of the model structure goes as
in \cite{model3}.  It is much easier since the category
$\dtop(\top_\Delta)$ is locally presentable by a proof similar to the
one of \cite[Proposition~6.11]{nonexistence}. Therefore, we do not
have to worry about the problems of smallness by
\cite[Proposition~1.3]{MR1780498}. So a big part of \cite{model3} can
be removed.  It remains to check that
\cite[Theorem~2.1.19]{MR99h:55031}:
\begin{itemize}
\item The class $\mathcal{W}$ of weak S-homotopy equivalences of
  $\dtop(\top_\Delta)$ satisfies the two-out-of-three property and is
  closed under retracts: clear.
\item $\cell(J^{gl}) \subset \mathcal{W} \cap \cof(I^{gl}_{+})$. The
  set inclusion $\cell(J^{gl}) \subset \cof(I^{gl}_{+})$ comes from
  the obvious set inclusion $J^{gl} \subset \cof(I^{gl})$.  The proof
  of the set inclusion $\cell(J^{gl}) \subset \mathcal{W}$ is similar
  to the proof of the same statement in \cite{model3}: the crucial
  facts are that:
\begin{itemize}
\item Proposition~\ref{calcul_explicite} which calculates a pushout by
  a map $\glob(U)\rightarrow \glob(V)$.
\item A map of the form $\id_X\p_\Delta j_n: X \p_\Delta \mathbf{D}^n
  \p_\Delta \{0\} \rightarrow X \p_\Delta \mathbf{D}^n \p_\Delta
  [0,1]$, where $j_n:\mathbf{D}^n \p_\Delta \{0\} \subset \mathbf{D}^n
  \p_\Delta [0,1]$ is the usual inclusion of spaces for some $n\geq
  0$, is a closed inclusion of a strong deformation retract of spaces.
  Indeed, since $\mathbf{D}^n$ is compact, the map $\id_A\p_\Delta
  j_n$ is equal to the map $\id_A\p_\mathcal{T} j_n$ by
  Proposition~\ref{produit}.  And we then consider once again the
  Str{\o}m model structure on the category of general topological
  spaces $\mathcal{T}$.
\item So for every map $f\in \cell(J^{gl})$, the continuous map $\P
  f:\P X \rightarrow \P Y$ is a closed inclusion of a strong
  deformation retract of spaces since the path space functor $\P:
  \dtop(\top_\Delta) \rightarrow \top_\Delta$ is finitely accessible
  as in Proposition~\ref{accessible}.  Therefore it is a weak homotopy
  equivalence.
\end{itemize}
\item $\inj(I^{gl}_{+})=\mathcal{W} \cap \inj(J^{gl})$. See
  \cite[Proposition~16.2]{model3} and \cite[Proposition~13.2]{model3}.
\end{itemize}
The model category $\dtop(\top_\Delta)$ is right proper since every
object is fibrant. The proof of the left properness of
$\dtop(\top_\Delta)$ is postponed until Appendix~\ref{leftmdtop}. It
is simplicial by the same proof as the one given in
\cite{realization}.  The second assertion of the statement is then
obvious. The functor $cat_\Delta:\mdtop_\Delta \rightarrow
\dtop(\top_\Delta)$ is colimit-preserving. Therefore it takes maps of
$\cell(I^{gl,top}_+)$ (resp.  $\cell(J^{gl,top})$) to maps of
$\cell(I^{gl}_+)$ (resp.  $\cell(J^{gl})$). Since every (resp.
trivial) cofibration of multipointed $d$-spaces is a retract of an
element of $\cell(I^{gl,top}_+)$ (resp. $\cell(J^{gl,top})$), one
deduces that the functor $cat_\Delta:\mdtop_\Delta \rightarrow
\dtop(\top_\Delta)$ takes (trivial) cofibrations to (trivial)
cofibrations. Let $f:X\rightarrow Y$ be a weak S-homotopy equivalence
between cofibrant multipointed $d$-spaces. The cofibrant object $X$ is
a retract of a globular complex $\overline{X}$. The space
$\P^{top}\overline{X}$ is compactly generated by
Proposition~\ref{bonspace1}. The space $\P \overline{X}$ is the path
space of the cofibrant flow $cat(\overline{X})$. Therefore the space
$\P \overline{X}$ is cofibrant by Proposition~\ref{path-cofibrant}.
Any cofibrant space is weak Hausdorff. So $\P \overline{X}$ is
compactly generated and by \cite[Theorem~IV.3.10]{model2}, the map
$\P^{top}\overline{X} \rightarrow \P \overline{X}$ is a weak homotopy
equivalence, and even a trivial Hurewicz fibration. Thus the map
$\P^{top}X \rightarrow \P X$ is a retract of the weak homotopy
equivalence $\P^{top}\overline{X} \rightarrow \P \overline{X}$. Hence
it is a weak homotopy equivalence as well. The commutative diagram of
$\Delta$-generated spaces
\[
\xymatrix{
\P^{top}X \fr{\simeq}\fd{\simeq} & \P^{top}Y \fd{\simeq}\\
\P X \fr{} & \P Y }
\] 
and the two-out-of-three property completes the proof. \epf 

In conclusion, let us say that the functor $cat_\Delta:\mdtop_\Delta
\rightarrow \dtop(\top_\Delta)$ is a kind of left Quillen equivalence
of cofibration categories.

\section{Underlying homotopy type of flows as a total left derived functor}
\label{under}

The underlying homotopy type functor of flows is defined in
\cite{model2}.  Morally speaking, it is the underlying topological
space of a flow, but it is unique only up to a weak homotopy
equivalence. It is equal, with the notations of this paper, to the
composite functor
\[\ho(\dtop(\top))\simeq \gltop[(\mathcal{W}^{top})^{-1}] \longrightarrow
\ho(\top_\Delta) \simeq \ho(\cgtop)\] where the middle functor is the
unique functor making the following diagram commutative:
\[
\xymatrix{
  \gltop \fr{X\mapsto |X|} \fd{} & \top_\Delta\fd{}\\
  \gltop/\!\sim_S\iso \gltop[(\mathcal{W}^{top})^{-1}] \fr{} &
  \ho(\top_\Delta)}
\] 
where both vertical maps are the canonical localization functors. 

\bp The functor $|-|:X=(|X|,X^0,\P^{top}X)\mapsto |X|$ from
$\mdtop_\Delta$ to $\top_\Delta$ is a left Quillen functor. \ep

\bpf By Proposition~\ref{limcolim}, this functor is topological. So it
has a right adjoint by \cite[Proposition~21.12]{topologicalcat}. In
fact, the right adjoint $R:\top_\Delta \rightarrow \mdtop_\Delta$ is
defined by:
\[R(Z)=(Z,\omega_\Delta(Z),\mtop_\Delta(([0,1],\{0,1\}),(Z,\omega_\Delta(Z)))).\]
The functor $|-|$ preserves cofibrations and trivial cofibrations by
Proposition~\ref{T1inc}.  \epf

\begin{cor} Up to equivalences of categories, the underlying homotopy
  type functor of flows is a total left derived functor.
\end{cor}

\bpf Indeed, the composite functor
\[\ho(\dtop(\top))\simeq \ho(\mdtop_\Delta) \stackrel{X\mapsto X^{cof}}\longrightarrow
\gltop/\!\sim_S \stackrel{Y\mapsto |Y|}\longrightarrow
\ho(\top_\Delta) \simeq \ho(\cgtop)\] is equal to the underlying
homotopy type functor.  \epf

The underlying space functor $X\mapsto |X|$ from $\mdtop_\Delta$ to
$\top_\Delta$ is not invariant with respect to weak S-homotopy
equivalences. With the identification $ \mathbf{S}^1=\{z\in
\mathbb{C},|z|=1\}$, consider the multipointed $d$-space
$X=(\mathbf{S}^1,\{1,\exp(i\pi/2)\}, \P^{top}X)$ where $\P^{top}X$ is
the closure by strictly increasing reparametrization of the set of
continuous paths $\{t\mapsto \exp(it\pi/2),t\in[0,1]\}$.  The
multipointed $d$-space has a unique initial state $1$ and a unique
final state $\exp(i\pi/2)$. Then consider the map of multipointed
$d$-spaces $f:\vI^{top} \longrightarrow X$ defined by
$|f|(t)=\exp(it\pi/2)$ for $t\in [0,1]$. Then $f$ is a weak S-homotopy
equivalence with $|\vI^{top}|$ contractible whereas $|X|$ is not so.

\appendix

\section{Left properness of $\dtop(\top_\Delta)$}
\label{leftmdtop}

As the proof of the left properness of $\dtop(\cgtop)$ given in
\cite[Section~7]{2eme}, the proof of the left properness of
$\dtop(\top_\Delta)$ lies in Proposition~\ref{calcul_explicite},
Proposition~\ref{versleft} and Proposition~\ref{petitplus}:

\bp \label{calcul_explicite} (Compare with
\cite[Proposition~7.1]{2eme}) Let $f:U\rightarrow V$ be a continuous
map between $\Delta$-generated spaces. Let $X$ be a flow enriched over
$\Delta$-generated spaces. Consider the pushout diagram of
multipointed $d$-spaces:
\[
\xymatrix{
\glob(U) \fd{\glob(f)}\fr{} & X \fd{} \\
\glob(V) \fr{} & \cocartesien Y.}
\] 
Then the continuous map $\P X \rightarrow \P Y$ is a transfinite
composition of pushouts of maps of the form
\[\id_{X_1} \p_\Delta \dots \p_\Delta f \p_\Delta \dots \p_\Delta \id_{X_p}\]
with $p\geq 0$.  \ep

\bpf The proof is exactly the same as for $\dtop(\cgtop)$. Nothing
particular need be assumed on the category of topological spaces we
are working with, except that it must be cartesian closed.  \epf

\bp \label{versleft} (Compare with \cite[Proposition~7.2]{2eme}) Let
$n\geq 0$. Let $i_n:\mathbf{S}^{n-1}\subset \mathbf{D}^n$ be the usual
inclusion of spaces. Let $X_1,\dots,X_p$ be $\Delta$-generated spaces.
Then the pushout of a weak homotopy equivalence along a map of the
form a finite product
\[\id_{X_1}\p_\Delta\dots\p_\Delta i_n\p_\Delta\dots \p_\Delta\id_{X_p}\]
with $p\geq 0$ is still a weak homotopy equivalence.  \ep

\bpf The three ingredients of \cite[Proposition~7.2]{2eme} are
\begin{enumerate}
\item The pushout of an inclusion of a strong deformation retract of
  spaces is an inclusion of a strong deformation retract of spaces.
  This assertion is true in the category of general topological spaces
  by \cite[Lemma~2.4.5]{MR99h:55031}.
\item The Seifert-Van-Kampen theorem for the fundamental groupoid
  functor which is true in the category of general topological spaces
  by \cite{0149.20002} \cite{MR2118984}.
\item The Mayer-Vietoris long exact sequence which holds in the
  category of general topological spaces: it is the point of view of,
  e.g., \cite{MR957919}.
\end{enumerate} \epf

The following notion is a weakening of the notion of closed
$T_1$-inclusion, introduced by Dugger and Isaksen.

\bd \cite[p 686]{MR2045835} An inclusion of spaces $f:Y\rightarrow Z$
is a {\rm relative $T_1$-inclusion of spaces} if for any open set $U$
of $Y$ and any point $z\in Z\backslash U$, there is an open set $W$ of
$Z$ with $U\subset W$ and $z\notin W$. \ed

\bp \label{i1} (Variant of \cite[Lemma~A.2]{MR2045835}) Consider a
pushout diagram of $\Delta$-generat\-ed spaces
\[
\xymatrix{
A\p_\Delta \mathbf{S}^{n-1} \fr{} \fd{} & Y \fd{} \\
A\p_\Delta \mathbf{D}^{n}  \fr{} & \cocartesien Z.}
\]
Then the map $Y\rightarrow Z$ is a relative $T_1$-inclusion of spaces.
\ep

\bpf By \cite[Lemma~A.2]{MR2045835}, the cocartesian diagram of
$\mathcal{T}$
\[
\xymatrix{
A\p_\mathcal{T} \mathbf{S}^{n-1} \fr{} \fd{} & Y \fd{} \\
A\p_\mathcal{T} \mathbf{D}^{n}  \fr{} & \cocartesien Z'.}
\]
yields a relative $T_1$-inclusion of spaces $Y\rightarrow Z'$. Since
both $\mathbf{S}^{n-1}$ and $\mathbf{D}^{n}$ are compact, one has
$A\p_\mathcal{T} \mathbf{S}^{n-1} \iso A\p_\Delta \mathbf{S}^{n-1}$
and $A\p_\mathcal{T} \mathbf{D}^{n} \iso A\p_\Delta \mathbf{D}^{n}$ by
Proposition~\ref{produit}.  Since colimits in $\mathcal{T}$ and in
$\top_\Delta$ are the same, the map $Y\rightarrow Z'$ is the map
$Y\rightarrow Z$.  \epf

\bp \cite[Lemma~A.3]{MR2045835} \label{finite} Any compact is finite
relative to the class of relative $T_1$-inclusions of spaces. \ep

\bp \label{petitplus} (Compare with \cite[Proposition~7.3]{2eme}) Let
$\lambda$ be an ordinal. Let $M:\lambda \longrightarrow \top_\Delta$
and $N:\lambda \longrightarrow \top_\Delta$ be two $\lambda$-sequences
of topological spaces.  Let $s:M\longrightarrow N$ be a morphism of
$\lambda$-sequences which is also an objectwise weak homotopy
equivalence. Finally, let us suppose that for all $\mu<\lambda$, the
continuous maps $M_\mu \longrightarrow M_{\mu+1}$ and $N_\mu
\longrightarrow N_{\mu+1}$ are pushouts of maps~\footnote{There is a
  typo error in the statement of \cite[Proposition~7.3]{2eme}. The
  expression ``pushouts of maps'' is missing.}  of the form of a
finite product \[\id_{X_1}\p_\Delta\dots\p_\Delta i_n\p_\Delta\dots
\p_\Delta\id_{X_p}\] with $p\geq 0$, with $i_n:\mathbf{S}^{n-1}\subset
\mathbf{D}^n$ being the usual inclusion of spaces for some $n\geq 0$.
Then the continuous map $\liminj s:\liminj M \longrightarrow \liminj
N$ is a weak homotopy equivalence.  \ep

\bpf The main ingredient of the proof of \cite[Proposition~7.3]{2eme}
is that any map $K \rightarrow M_\lambda$ factors as a composite
$K\rightarrow M_\mu \rightarrow M_\lambda$ for some ordinal
$\mu<\lambda$ as soon as $K$ is compact if $\lambda$ is a limit
ordinal.  More precisely, one needs to apply this fact for $K$
belonging to the set $\{\mathbf{S}^n,\mathbf{S}^n\p_\Delta [0,1],n\geq
0\}$. It is not true that the maps $M_\mu \rightarrow M_{\mu+1}$ and
$N_\mu \rightarrow N_{\mu+1}$ are closed $T_1$-inclusions of spaces
since the spaces $X_i$ are not necessarily weak Hausdorff anymore. So
we cannot apply \cite[Proposition~2.4.2]{MR99h:55031} saying that
compact spaces are finite relative to closed $T_1$-inclusions of
spaces, unlike in the proof of \cite[Proposition~7.3]{2eme}.  Let
$Z=X_1 \p_\Delta \dots \p_\Delta X_p$. Then each map $M_\mu
\longrightarrow M_{\mu+1}$ and $N_\mu \longrightarrow N_{\mu+1}$ is a
pushout of a map of the form $Z \p_\Delta \mathbf{S}^{n-1} \rightarrow
Z \p_\Delta \mathbf{D}^{n}$. We can then apply Proposition~\ref{i1}
and Proposition~\ref{finite}. The proof is therefore complete. \epf

\section{The simplicial structure of $\mdtop_\Delta$}
\label{simp}

The construction is very similar to the one given in
\cite{realization} for the category of flows.

\bd Let $K$ be a non-empty connected simplicial set. Let $X$ be an
object of $\mdtop_\Delta$. Let $X \ot K=|K| \boxtimes X$ where $|K|$
means the geometric realization of $K$ \cite{MR2001d:55012}.  \ed

\bd Let $K$ be a non-empty simplicial set. Let $(K_i)_{i\in I}$ be its
set of non-empty connected components. Let $X \ot K := \bigsqcup_{i\in
  I} X \ot K_i$.  And let $X \ot \varnothing = \varnothing$.  \ed

\bp Let $K$ be a simplicial set. Then the functor $-\ot K:\mdtop_\Delta
\rightarrow \mdtop_\Delta$ is a left adjoint. \ep

\bpf As in \cite{realization}, it suffices to prove the existence of
the right adjoint \[(-)^ K:\mdtop_\Delta \rightarrow \mdtop_\Delta\]
for $K$ non-empty connected and to set:
\begin{itemize}
\item $X ^\varnothing = \mathbf{1}$
\item for a general simplicial set $K$ with non-empty connected
  components $(K_i)_{i\in I}$, let $X^K := \prod_{i\in I} X^{K_i}$.
\end{itemize}
So now suppose that $K$ is non-empty connected. For a given
multipointed $d$-space $X$, let (compare with
\cite[Notation~7.6]{model3} and \cite[Theorem~7.7]{model3}):
\begin{itemize}
\item $(X^K)^0=X^0$
\item $|X^K|=\ttop_\Delta(|K|,|X|)$
\item for $(\alpha,\beta)\in X^0\p X^0$,
  $\P_{\alpha,\beta}^{top}(X^K)=\ttop_\Delta(|K|,\P_{\alpha,\beta}^{top}X)$.
\end{itemize}
We can observe that the functor $(-)^ K:\mdtop_\Delta \rightarrow
\mdtop_\Delta$ commutes with limits and is $\lambda$-accessible if
$|K|$ is $\lambda$-presentable in $\top_\Delta$ for some regular
cardinal $\lambda$ since the functor $\P^{top}:\mdtop_\Delta
\rightarrow \set$ is finitely accessible by
Proposition~\ref{accessible}. So by \cite[Theorem~1.66]{MR95j:18001},
it is a right adjoint. It is easy to check that the left adjoint is
precisely $-\ot K$.  \epf

\bp \label{map} Let $X$ and $Y$ be two multipointed $d$-spaces. Let
$\Delta[n]$ be the $n$-simplex. Then there is a natural isomorphism of
simplicial sets \[\mdtop_\Delta(X\ot \Delta[*],Y) \iso \sing
\mdttop_\Delta(X,Y)\] where $\sing$ is the singular nerve functor.
This simplicial set is denoted by $\map(X,Y)$.  \ep

\bpf Since $\Delta[n]$ is non-empty and connected, one has
\[\sing(\mdttop_\Delta(X,Y))_n=\top_\Delta(|\Delta[n]|,\mdttop_\Delta(X,Y))
\iso \mdtop_\Delta(X\ot \Delta[n],Y).\] \epf

\bth\label{simpl} The model category $\mdtop_\Delta$ together
with the functors $-\ot K$, $(-)^K$ and $\map(-,-)$ assembles to a
simplicial model category.  \eth

\bpf Proof analogous to the proof of
\cite[Theorem~3.3.15]{realization}.  \epf

In fact, the category of multipointed $d$-spaces $\mdtop_\Delta$ is
also tensored and cotensored over $\top_\Delta$ in the sense of
\cite{strom2} because of Proposition~\ref{connexe} and
Theorem~\ref{cyl}. On the contrary, one has:

\bp The category $\mdtop_\mathcal{K}$ of multipointed $d$-spaces over
$\top_\mathcal{K}$ is neither tensored, nor cotensored over
$\top_\mathcal{K}$. \ep

\bpf Otherwise the functor $\mdttop_\mathcal{K}(X,-) :
\mdtop_\mathcal{K} \rightarrow \top_\mathcal{K}$ would preserve
limits. Take $X=\{0\}$. Then for any multipointed $d$-space $Y$,
$\mdttop_\mathcal{K}(\{0\},Y)$ is the discrete space $Y^0$ by
Proposition~\ref{map}. But a limit of discrete spaces in
$\top_\mathcal{K}$ is not necessarily discrete (e.g. the $p$-adic
integers $\Z_p =\limproj \Z/p^n\Z$ \cite{topologie}). Contradiction.
\epf

The same phenomenon arises for the category of flows: read the comment
\cite[p567]{model3} after the statement of Theorem~5.10.

\bibliographystyle{alpha}
\bibliography{Mdtop}

\end{document}